\documentclass[12pt]{article}
\author{Bing-Long Chen and Xi-Ping Zhu\\[8pt]
Department of Mathematics \\
Zhongshan University \\
Guangzhou, P.R.China}
\title{\textbf{Uniqueness of the Ricci Flow on Complete Noncompact Manifolds}}
\date {May, 2005}
\usepackage{latexsym}
\usepackage{amsmath}
\usepackage{amssymb,amscd}
\usepackage[mathscr]{eucal}
\newtheorem{thm}{Theorem}[section]
\newtheorem{cor}[thm]{Corollary}
\newtheorem{lem}[thm]{Lemma}
\newtheorem{prop}[thm]{Proposition}

\numberwithin{equation}{section}
\begin{document}
\maketitle

\centerline{\Large{\textbf{Abstract}}}
 \vskip 0.5cm

  The Ricci flow is an evolution system on metrics. For a given
metric as initial data, its local existence and uniqueness on
compact manifolds was first established by Hamilton \cite{Ha1}.
Later on, De Turck \cite{De} gave a simplified proof. In the later
of 80's, Shi \cite{Sh1} generalized the local existence result to
complete noncompact manifolds. However, the uniqueness of the
solutions to the Ricci flow on complete noncompact manifolds is
still an open question. Recently it was found that the uniqueness
of the Ricci flow on complete noncompact manifolds is important in
the theory of the Ricci flow with surgery. In this paper, we give
an affirmative answer for the uniqueness question. More precisely,
we prove that the solution of the Ricci flow with bounded
curvature on a complete noncompact manifold is unique.

\pagebreak[4]
  \section{Introduction}
\vskip 0.5cm

\qquad  Let $(M^n, g_{ij})$ be a complete Riemannian (compact or
noncompact) manifold. The Ricci flow
$$
\frac{\partial}{\partial t}g_{ij}(x,t)=-2R_{ij}(x,t), \ \ \ \ \ \
\ \ \ \text{for}\ x\in M^{n} \text{ and }\  t\geq 0, \eqno{(1.1)}
$$
with $g_{ij}(x,0)=g_{ij}(x)$, is a weakly parabolic system on
metrics. This evolution system was introduced by Hamilton in
\cite{Ha1}. Now it has proved to be powerful in the research of
differential geometry and lower dimensional topology (see for
example Hamilton's works \cite{Ha1}, \cite{Ha2}, \cite{Ha3},
\cite{Ha6} and the recent works of Perelman \cite{P1}, \cite{P2}).
The first matter for the Ricci flow (1.1) is the short time
existence and uniqueness of the solutions. When the manifold
$M^{n}$ is compact, Hamilton proved in \cite{Ha1} that the Ricci
flow (1.1) has a unique solution for a short time. So the problem
has been well settled on compact manifolds. In \cite{De}, De Turck
introduced an elegant trick to give a simplified proof. Later on,
Shi \cite{Sh1} extended the short time existence result to
noncompact manifolds. More precisely, Shi \cite{Sh1} proved that
if $(M^n, g_{ij})$ is complete noncompact with bounded curvature,
then the Ricci flow (1.1) has a solution with bounded curvature on
a short time interval. In this paper, we will deal with the
uniqueness of the Ricci flow on complete noncompact manifolds.

The uniqueness of the Ricci flow is important in the theory of the
Ricci flow with surgery (see for example \cite{P1}, \cite{P2} and
\cite{CZ1}). When we consider the Ricci flow on a compact
manifold, the Ricci flow will generally develop singularities in
finite time. In the theory of the Ricci flow with surgery, one
eliminates the singularities by Hamilton's geometric surgeries
(cut off the high curvature part and glue back a standard cap,
then run the Ricci flow again). An important question in this
theory is to control the curvature of the glued cap after surgery.
The uniqueness theorem of the Ricci flow insures that the solution
on glued cap is sufficiently close to a (complete noncompact)
standard solution, which is the evolution of capped round
cylinder. Then we can apply the estimate of the standard solutions
\cite{P2} and \cite{CZ1} to get the desired control on curvature.
The employing of the uniqueness theorem is essential. So even if
we consider the Ricci flow on compact manifolds, we still have to
encounter the problem of uniqueness on noncompact manifolds.

It is well-known that the uniqueness of the solution of a
parabolic system on a complete noncompact manifold does not always
hold if one does not impose any growth condition of the solutions.
For example, even the simplest linear heat equation on
$\mathbb{R}$ with zero as initial data has a nontrivial solution
which grows faster than $e^{a|x|^{2}}$ for any $a>0$ whenever
$t>0$. This says, for the standard linear heat equation, the most
growth rate for the uniqueness is $e^{a|x|^{2}}$. Note that in a
K$\ddot{a}$hler manifold, the Ricci curvature is given by
$$R_{i\bar{j}} = -\frac{\partial^2}{\partial z^i \partial
\bar{z}^j}log det (g_{k\bar{l}}).$$ Thus the reasonable growth
rate that we can expect for the uniqueness of the Ricci flow is
the solution with bounded curvature.

In this paper, we will prove the following uniqueness theorem of
the Ricci flow.
\begin{thm} Let $(M^n,g_{ij}(x))$ be a complete noncompact Riemannian manifold
of dimension $n$ with bounded curvature. Let ${g}_{ij}(x,t)$ and
$\bar{g}_{ij}(x,t)$ be two solutions to the Ricci flow on
$M^n\times[0,T]$ with the same $g_{ij}(x)$ as initial data and
with bounded curvatures. Then $g_{ij}(x,t)= \bar{g}_{ij}(x,t)$ for
all $(x,t)\in M^{n}\times [0,T]$.
\end{thm}

 Since the Ricci flow is not a strictly parabolic system, our argument
 will apply the De Turck trick. This is to consider the composition
  of the Ricci flow with a family of diffeomorphisms generated by
  the harmonic map flow. By pulling back the Ricci flow by this family of diffeomorphisms,
   the evolution equations become
  strictly parabolic. In order to use the uniqueness theorem of
   a strict parabolic system on a noncompact manifold, we have to overcome
   two difficulties. The first one is to establish a short time
   existence for the harmonic map flow between noncompact
   manifolds. The second one is to get a priori estimates for the
   harmonic map flow so that after pulling backs, the solutions to
   the strictly parabolic system still satisfy suitable growth
   conditions. To the best of our knowledge, one can only get short time
   existence of harmonic map flow by imposing negative curvature
   or convex condition on the target manifolds (see for example, \cite{ES} and \cite{DL}) or by simply
   assuming the image of initial data lying in a compact domain on
   the target manifold (see for example \cite{LT}). In \cite{CZ1}, we observed that the condition of
   injectivity radius bounded from below ensures certain uniform
   (local) convexity and showed that this is sufficient to give
   the short time existence and the a priori estimates for the
   harmonic map flow. Thus in \cite{CZ1}, we obtained the uniqueness under
   an additional assumption that the initial metric has a positive
   lower bound on injectivity radius. The main purpose of this
   paper is to remove this additional assumption. Note
   from \cite{CLY} or \cite{CGT} that the injectivity radius of the initial manifold
   decays at worst exponentially. This allows us to conformally
   straighten the initial manifold at infinity. Our idea is to
   study the evolution equations coming from the composition of
   the Ricci flow and harmonic map flow, as well as a conformal
   change.

   This new approach has the advantage of transforming the Ricci
   flow equation to a strictly parabolic system
   on a manifold with uniform geometry at infinity. We expect that it
   could also give new short time existence for the Ricci flow without
   assuming the boundedness of the curvature of the initial
   metric. As a direct consequence, we have the following result.
   \begin{cor}
   Suppose $(M^{n},g_{ij}(x))$ is a complete Riemannian manifold,
   and suppose $g_{ij}(x,t)$ is a solution to the Ricci flow with bounded
   curvature on $M^{n}\times[0,T]$ and with $g_{ij}(x)$ as initial
   data. If $G$ is the isometry group of $(M^{n},g_{ij}(x))$, then
   $G$ remains to be an isometric subgroup of
   $(M^{n},g_{ij}(x,t))$ for each $t\in[0,T]$.
   \end{cor}

   This paper is organized as follows. In Section 2, we study the
   harmonic map flow coupled with the Ricci flow. In Section 3, we
   study the Ricci-De Turck flow and prove the uniqueness theorem.

   We are grateful to Professor S. T. Yau
for many helpful discussions and encouragement. The second author
is partially supported by the IMS of The Chinese University of
Hong Kong and the first author is supported by FANEDD 200216 and
NSFC 10401042.

\vskip 0.8cm
\section{Harmonic map flow coupled with the Ricci flow}

\vskip 0.5cm
    Let $(M^{n},g_{ij}(x))$ and $(N^{m},h_{ij}(y))$ be two
    Riemannian manifolds,
     $f:M^{n}\rightarrow N^{m}$ be a map. The harmonic map
     flow is the following evolution equation for maps from $M^n$
     to $N^m$,
\begin{equation*}\tag{2.1}
\left\{
\begin{split}
 \quad \frac{\partial}{\partial t}f(x,t)&=\triangle f(x,t), \ \ \ \mbox{ for }x \in M^n, t>0,  \\
  f(x,0)&= f(x), \  \  \  \  \   \mbox{ for }x \in M^n ,
  \end{split}
 \right.
\end{equation*}
 where $\triangle$ is defined by using the metric $g_{ij}(x)$ and
$h_{\alpha\beta}(y)$ as follows
$$
 \triangle
f^{\alpha}(x,t)=g^{ij}(x)\nabla_{i}\nabla_{j}f^{\alpha}(x,t),
   \ \ \ \ \ \ \ \ \ \ $$
   and
   $$
   \nabla_{i}\nabla_{j}f^{\alpha}=\frac{\partial^{2}f^{\alpha}}{\partial
    x^{i}\partial
   x^{j}}-\Gamma^{k}_{ij}\frac{\partial f^{\alpha}}{\partial
   x^{k}}+\Gamma^{\alpha}_{\beta\gamma}\frac{\partial f^{\beta}}{\partial x^{i}}
   \frac{\partial f^{\gamma}}{\partial x^{j}}. \eqno(2.2)$$
   Here we use $\{x^{i}\}$ and $\{y^{\alpha}\}$ to denote the
   local
   coordinates of $M^{n}$ and $N^{m}$
   respectively, $\Gamma^{k}_{ij}$ and
   $\Gamma^{\alpha}_{\beta\gamma}$ the corresponding Christoffel
   symbols of $g_{ij}$ and $h_{\alpha\beta}$.

Let $g_{ij}(x,t)$ be a complete smooth solutions of the Ricci flow
with $g_{ij}(x)$ as initial data, then the harmonic map flow
coupled with Ricci flow is the following equation
\begin{equation*}
\left\{
\begin{split}
 \quad \frac{\partial}{\partial t}f(x,t)&=\triangle_{t} f(x,t), \ \ \ \mbox{ for }x \in M^n, t>0,  \\
  f(x,0)&= f(x), \  \  \  \  \   \mbox{ for }x \in M^n ,
  \end{split}
 \right.
\end{equation*}
where $\triangle_{t}$ is defined as above by using the metric
$g_{ij}(x,t)$ and $h_{\alpha\beta}(y)$.

Suppose $g_{ij}(x,t)$ is a solution to the Ricci flow on $
M^{n}\times [0,T]$ with bounded curvature $$ |Rm|(x,t)\leq k_{0}$$
 for all $(x,t)\in M^{n}\times[0,T]$. Let $(N^{n},h_{\alpha\beta})=(M^{n},g_{ij}(\cdot,T))$
 be the target manifold. The purpose of this section is to prove
 the following theorem
     \begin{thm}
     There exists $0<T_0<T$, depending only on $k_{0}$, $T$ and $n$ such that the harmonic map flow
   coupled with the Ricci flow
     \begin{equation*}\tag{2.3}
     \left\{
     \begin{split}
     \quad \frac{\partial}{\partial t}F(x,t)&=\triangle_{t} F(x,t),  \\
  F(\cdot,0)&= identity, \  \  \  \  \
  \end{split}
  \right.
     \end{equation*}
     has a solution on $M^{n}\times[0,T_0]$ satisfying the
     following estimates
     $$
     \begin{array}{lr}
     |\nabla F|\leq \tilde{C}_{1},\\
           |\nabla^{k} F|\leq \tilde{C}_{k} t^{-\frac{k-2}{2}},\ \ \ \ \
      \mbox{for all k }\ \geq 2,
     \end{array}\eqno{(2.4)}
     $$
     for some constants $\tilde{C}_{k}$ depending only on $k_{0}$,
     $T$, $k$ and $n$.
     \end{thm}

     The proof will occupy the rest of this section.

\subsection{Expanding base and target metrics at infinity}

We will construct appropriate auxiliary functions on $M^{n}$ and
$N^{n}$ and do conformal deformations for the base and the target
metrics. Firstly, we construct the function on
$(N^{n},h_{\alpha\beta})$. The function can be obtained by solve
certain equations \cite{ScY}
 or smoothing certain functions by convolution \cite{GW}.
\begin{lem}
Fix $p\in N^{n}$. Then for any $a\geqslant 1$, there exists a
$C^{\infty}$ nonnegative function $\varphi_{a}$ on $N^{n}$ such
that
\begin{equation*}\tag{2.5}
     \left\{
     \begin{split}
      \varphi_{a}(y)&\equiv 0  \ \ \  &&\text{on}\ \ \ \  B(p,a),  \\
  d(y,p)&\leqslant\varphi_{a}(y)\leqslant C_{0}d(y,p) \ \ \  &&\text{on}\ \ \ \ N^{n}\backslash B(p,2a),\\
 |\nabla^{k}\varphi_{a}| &\leqslant C_{k} \ \ \  && \text{on}\ \ \ \
 N^{n},\ \ \ \ \text {for}\  k\geqslant1,
  \end{split}
  \right.
     \end{equation*}
     where $C_{k}$, $i=0, 1, 2, \cdots,$ are constants depending only on $k_{0}$ and $T$; the distance $d(y,p)$, the covariant
     derivatives $\nabla^{k}\varphi_{a}$ and the norms
     $|\nabla^{k}\varphi_{a}|$ are computed by using the metric
     $h_{\alpha\beta}$.
\end{lem}
Proof. Let $\xi$ be a smooth nonnegative increasing function on
$\mathbb{R}$ such that $\xi(s)=0$ for $s\in
(-\infty,\frac{5}{4}]$, and $\xi=1$ for
$s\in[\frac{7}{4},\infty)$. For each $y\in N^{n}$, by averaging
the functions $\xi(\frac{d(p,y)}{a})$ and $d(p,y)$ over a suitable
ball of the tangent space $T_{y}N^{n}$ (see for example
\cite{GW}), we obtain two smooth functions $\xi_{a}$ and $\rho$.
Notice that $(N^{n},h_{\alpha\beta})=(M^{n},g_{ij}(\cdot,T))$,
thus all the covariant derivatives of the curvatures of
$h_{\alpha\beta}$ are bounded by using Shi's gradient estimates
\cite{Sh1}. Then $\varphi_{a}=C\xi_{a}\rho $, for some constant
$C$ depending only on $k_{0}$ and $T$, is the desired function.

 $\hfill\#$

Recall from \cite{CLY} and \cite{CGT} that on a complete manifold
with bounded curvature, the injectivity radius decays at worst
exponentially; more precisely, there exists a constant
$\tilde{C}(n)>0$ depending only on the dimension, and there exists
a constant $\delta>0$ depending on $n$, $k_{0}$ and the
injectivity radius at $p$ such that
$$
inj(N^{n},h_{\alpha\beta},y)\geqslant \delta
e^{-\tilde{C}(n)\sqrt{k_{0}}d(y,p)}. \eqno(2.6)
$$

Fix $a\geqslant1$, let
$\varphi^{a}=4\tilde{C}(n)\sqrt{k_{0}}\varphi_{a}$ and set
$$
h^{a}_{\alpha\beta}=e^{\varphi^{a}}h_{\alpha\beta}. \eqno(2.7)
$$
Clearly, $h^{a}_{\alpha\beta}=h_{\alpha\beta}$ on $B(p,a)$. Note
that $(N^{n},h_{\alpha\beta})=(M^{n},g_{ij}(\cdot,T))$, so the
function $\varphi_{a}$ is also a function on $M^{n}$. Let
$$
g^{a}_{ij}(x,t)=e^{\varphi^{a}}g_{ij}(x,t) \eqno(2.8)
$$
be the new family of metrics on $M^{n}$. Instead of (2.3), we will
consider a new harmonic map flow
$$ \left\{
     \begin{split}
     \quad \frac{\partial}{\partial t}\overset{a}{F}(x,t)&=\overset{a}{\triangle}_{t} \overset{a}{F}(x,t),  \\
  \overset{a}{F}(\cdot,0)&= identity, \  \  \  \  \
  \end{split}
  \right.\eqno(2.3)_{a}
  $$
where $\overset{a}{\triangle}_{t}\overset{a}{F}$ is defined  by
using the metric $g^{a}_{ij}(x,t)$ and $h^{a}_{\alpha\beta}(y)$.

Before we can solve $(2.3)_{a}$, we have to discuss the geometry
of the new metrics  $h^{a}_{\alpha\beta}(y)$ and
$g^{a}_{ij}(x,t)$. Let us first compute the curvature and its
covariant derivatives and injectivity radius of
$(N^{n},h^a_{\alpha\beta})$ as follows.

By a direct computation, we get
\begin{equation*}\tag{2.9}
\begin{split}
\overset{a}{R}_{\alpha\beta\gamma\delta}=&e^{\varphi^{a}}
R_{\alpha\beta\gamma\delta}+\frac{e^{\varphi^{a}}}{4}\{|\nabla
\varphi^{a}|^{2}(h_{\alpha\delta}h_{\beta\gamma}-h_{\alpha\gamma}
h_{\beta\delta})\\&+(2\nabla_{\alpha}\nabla_{\delta}\varphi^{a}-
\nabla_{\alpha}\varphi^{a}\nabla_{\delta}\varphi^{a})h_{\beta\gamma}
+(2\nabla_{\beta}\nabla_{\gamma}\varphi^{a}-
\nabla_{\beta}\varphi^{a}\nabla_{\gamma}\varphi^{a})h_{\alpha\delta}\\&-(2\nabla_{\beta}\nabla_{\delta}\varphi^{a}-
\nabla_{\beta}\varphi^{a}\nabla_{\delta}\varphi^{a})h_{\alpha\gamma}-(2\nabla_{\alpha}\nabla_{\gamma}\varphi^{a}-
\nabla_{\alpha}\varphi^{a}\nabla_{\gamma}\varphi^{a})h_{\beta\delta}\}
\end{split}
\end{equation*}
where $\overset{a}{R}_{\alpha\beta\gamma\delta}$ is the curvature
of $h^{a}_{\alpha\beta}$,
$\nabla_{\alpha}\varphi^{a}$,$\nabla_{\alpha}\nabla_{\delta}\varphi^{a}$
and $|\nabla_{\alpha}\varphi^{a}|$ are computed by the metric
$h_{\alpha\beta}$. Therefore, by combining with (2.5), we have
\begin{equation*}\tag{2.10}
\begin{split}
|\overset{a}{R}_{m}|_{h^{a}}&\leqslant
e^{-\varphi^{a}}(k_{0}+C(n)(C_{2}+C_{1}^{2}))\\
&<\infty.
\end{split}
\end{equation*}

For higher derivatives, we rewrite (2.9) in a simple form
$$
\overset{a}{R}_{m}=e^{\varphi^{a}}\{R_{m}+\nabla\varphi^{a}\ast\nabla\varphi^{a}\ast
h^{2}\ast h^{-1}+\nabla^{2}\varphi^{a}\ast h\}
$$
where we use $A\ast B$ to express some linear combinations of
tensors formed by contractions of tensor product of $A$ and $B$.
Note that
\begin{equation*}
\begin{split}
\overset{a}{\Gamma^{\alpha}_{\beta\gamma}}-\Gamma^{\alpha}_{\beta\gamma}
&=\frac{1}{2}[\nabla_{\beta}\varphi^{a}\delta^{\alpha}_{\gamma}+\nabla_{\gamma}\varphi^{a}
\delta^{\alpha}_{\beta}-h^{\alpha\eta}h_{\beta\gamma}\nabla_{\eta}\varphi^{a}]\\
&=(\nabla\varphi^{a}\ast h\ast h^{-1})^{\alpha}_{\beta\gamma},
\end{split}
\end{equation*}
so by induction, we have
\begin{equation*}\tag{2.11}
\begin{split}
\overset{a}{\nabla^{k}}\overset{a}{R}_{m}&=\nabla\overset{a}{{\nabla}^{k-1}}
\overset{a}{R_{m}}+(\overset{a}{\Gamma}-\Gamma)\ast\overset{a}{\nabla^{k-1}}\overset{a}{R_{m}}\\
&=e^{\varphi^{a}}\{\sum^{k}_{l=0}\nabla^{l}R_{m}\ast\sum_{i_{1}+\cdots+i_{p}=k-l}\nabla^{i_{1}}\varphi^{a}\ast
\cdots\ast
\nabla^{i_{p}}\varphi^{a}+\sum_{i_{1}+\cdots+i_{p}=k+2}\nabla^{i_{1}}\varphi^{a}\ast
\cdots\ast \nabla^{i_{p}}\varphi^{a}\},
\end{split}
\end{equation*}
where we denote $\nabla^{0}\varphi^{a}=1$. By combining with (2.5)
and gradient estimate of Shi \cite{Sh1}, we get
\begin{equation*}\tag{2.12}
\begin{split}
|\overset{a}{\nabla^{k}}\overset{a}{R}_{m}|_{h^{a}}&\leqslant
e^{-\frac{k+2}{2}\varphi^{a}}C(n,k_{0},k,C_{1},\cdots,C_{k+2})(\sum^{k}_{l=0}|\nabla^{l}R_{m}|+1)\\
&\leqslant C(n,k_{0},T,k)e^{-\frac{k+2}{2}\varphi^{a}}\\
&\leqslant C(n,k_{0},T,k).
\end{split}
\end{equation*}

For the injectivity radius of $h^{a}_{\alpha\beta}$, we know from
(2.5) and (2.7) that for any $y\in N^{n}\backslash B(p,2a+1)$,
$$
\overset{a}{B}(y,1)\supset
B(y,e^{-2\tilde{C}(n)\sqrt{k_{0}}(\varphi_{a}+C_{1})})
$$
and
\begin{equation*}\tag{2.13}
\begin{split}
Vol_{h^{a}}(\overset{a}{B}(y,1))&=\int_{\overset{a}{B}(y,1)}(e^{4\tilde{C}(n)\sqrt{k_{0}}\varphi_{a}})^{\frac{n}{2}}
\\
&\geqslant
e^{2\tilde{C}(n)\sqrt{k_{0}}(\varphi_{a}-C_{1})}Vol_{h}(B(y,e^{-2\tilde{C}(n)\sqrt{k_{0}}(\varphi_{a}+C_{1})}))
\end{split}
\end{equation*}
where we denote by $\overset{a}{B}(y,1)$ the ball centered at $y$
and of radius 1 with respect to metric $h^{a}_{\alpha\beta}$, and
$Vol_{h^{a}}(\overset{a}{B}(y,1))$ its volume.

Since $$\varphi_{a}(y)\geqslant d(y,p),$$ for $y\in
N^{n}\backslash B(p,2a+1)$, there holds
$$
e^{-2\tilde{C}(n)\sqrt{k_{0}}(\varphi_{a}+C_{1})}\leqslant \delta
e^{-\tilde{C}(n)\sqrt{k_{0}}d(y,p)} \eqno(2.14)
$$
for $y\in N^{n}\backslash
B(p,2a+1+|\frac{log\delta^{-1}}{\tilde{C}(n)\sqrt{k_{0}}}|)$. By
(2.6), (2.10), (2.13), (2.14) and volume comparison theorem, we
have
\begin{equation*}\tag{2.13}
\begin{split}
Vol_{h^{a}}(\overset{a}{B}(y,1))&\geqslant c(n,k_{0})
e^{2n\tilde{C}(n)\sqrt{k_{0}}(\varphi_{a}-C_{1})}(e^{-2\tilde{C}(n)\sqrt{k_{0}}(\varphi_{a}+C_{1})})^{n}
\\
&\geqslant c(n,k_{0}),
\end{split}
\end{equation*}
By combining this with the local injectivity radius estimate in
\cite{CLY} or \cite{CGT}, we get
$$
inj(N^{n},h^{a},y)\geqslant \tilde{C}(n,k_{0})>0, \ \ \ \mbox{for}
\ \ y\in N^{n}\backslash
B(p,2a+1+|\frac{log\delta^{-1}}{\tilde{C}(n)\sqrt{k_{0}}}|).
$$
Consequently, we have proved the  following lemma.
\begin{lem}
There exists a sequence of constants $\bar{C_{0}}$, $
\bar{C_{1}}$, $ \cdots$, with the following property. For all
$a\geqslant 1$, there exists $i_{a}>0$, such that the metrics
$h^{a}_{\alpha\beta}=e^{\varphi^{a}}h_{\alpha\beta}$ on $N^{n}$
satisfy
\begin{equation*}\tag{2.15}
\begin{split}
|\overset{a}{\nabla^{k}}\overset{a}{R}_{m}|_{h^{a}}&\leqslant \bar{C}_{k}e^{-\frac{k+2}{2}\varphi^{a}}\leqslant \bar{C}_{k}\\
inj(N^{n},h^{a}_{\alpha\beta})&\geqslant i_{a}>0
\end{split}
\end{equation*}
for $k=0, 1, \cdots.$
\end{lem}

 $\hfill\#$

We next estimate the curvature and the its covariant derivatives
of $g^{a}_{ij}(x,t)=e^{\varphi^{a}}g_{ij}(x,t)$.

By the Ricci flow equation, we have
\begin{equation*}\tag{2.16}
\begin{split}
\Gamma^{l}_{ij}(\cdot,T)-\Gamma^{l}_{ij}(\cdot,t)&=\int_{t}^{T}(g^{-1}\ast\nabla
Ric)(\cdot,s)ds,\\
  \nabla^{k}_{g(\cdot,T)}(\Gamma(\cdot,T)-\Gamma(\cdot,t))&=\int_{t}^{T}\sum^{k}_{l=0}\nabla^{k+1-l}Ric
 \ast
 \sum_{i_1+1+\cdots+i_{p}+1=l}\nabla^{i_1}_{g(\cdot,T)}(\Gamma(\cdot,s)-\Gamma(\cdot,T))\\&
 \ast\cdots\ast\nabla^{i_p}_{g(\cdot,T)}(\Gamma(\cdot,s)-\Gamma(\cdot,T))\ast
 g^{k}\ast g^{-(k+1)}(\cdot,s)ds.
\end{split}
\end{equation*}
By combining with the gradient estimates of Shi \cite{Sh1} and
induction on $k$, we have
\begin{equation*}\tag{2.17}
\left\{
\begin{split}
|\Gamma(\cdot,T)-\Gamma(\cdot,t)|&\leqslant
C(n,k_{0},T)\int_{t}^{T}\frac{1}{\sqrt{s}}ds,\\
  |\nabla_{g(\cdot,T)}(\Gamma(\cdot,T)-\Gamma(\cdot,t))|&\leqslant
  C(n,k_{0},T)(1+|log t|),\\
|\nabla^{k}_{g(\cdot,T)}(\Gamma(\cdot,T)-\Gamma(\cdot,t))|&\leqslant
  C(n,k_{0},T,k)t^{-\frac{k-1}{2}}, \ \ \ \mbox{ for } k\geq2.\\
\end{split}
\right.
\end{equation*}
Since
$$
\nabla^{k}_{g(\cdot,t)}\varphi^{a}=\sum^{k-1}_{l=0}\nabla_{g(\cdot,T)}^{k-l}\varphi^{a}\ast
\sum_{i_1+1+\cdots+i_p+1=l}\nabla^{i_1}_{g(\cdot,T)}(\Gamma(\cdot,t)-\Gamma(\cdot,T))
 \ast\cdots\ast\nabla^{i_p}_{g(\cdot,T)}(\Gamma(\cdot,t)-\Gamma(\cdot,T))
$$
for $k\geqslant1$, the combination with (2.17) and (2.5) gives
\begin{equation*}\tag{2.18}
\left\{
\begin{split}
&|\nabla_{g(\cdot,t)}\varphi^{a}|+|\nabla_{g(\cdot,t)}^{2}\varphi^{a}|\leqslant
C(n,k_{0},T),\\
  &|\nabla_{g(\cdot,t)}^{3}\varphi^{a}|\leqslant
  C(n,k_{0},T)(1+|log t|),\\
&|\nabla^{k}_{g(\cdot,t)}\varphi^{a}|\leqslant
  C(n,k_{0},T,k)t^{-\frac{k-3}{2}}, \ \ \ \ \text{for}\ \  k\geqslant4.\\
\end{split}
\right.
\end{equation*}
Then by combining (2.11) and (2.18), the curvature and the
covariant derivatives of $g^{a}(\cdot,t)$ can be estimated as
follows
\begin{equation*}
\begin{split}
|\overset{a}{\nabla}^{k}\overset{a}{R}_{m}|_{g^{a}(\cdot,t)}\leqslant
  C(n,k_{0},T,k) e^{-\frac{k+2}{2}\varphi^{a}}t^{-\frac{k}{2}}, \ \ \ \ \text{for}\ \  k\geqslant 0.\\
\end{split}
\end{equation*}
Summing up, the above estimates give the following
\begin{lem}
There exists a sequence of constants $\bar{k_{0}}$, $
\bar{k_{1}}$,  $ \cdots$, with the following property. For all
$a\geqslant 1$, the metrics
$g^{a}_{ij}(\cdot,t)=e^{\varphi^{a}}g_{ij}(\cdot,t)$ on $M^{n}$
satisfy
\begin{equation*}\tag{2.19}
\begin{split}
|\overset{a}{\nabla^{l}}\overset{a}{R}_{m}|_{g^{a}(\cdot,t)}
&\leqslant \bar{k}_{l}
e^{-\frac{k+2}{2}\varphi^{a}}t^{-\frac{l}{2}}, \ \ \ \text{for}\ \
\ l\geqslant 0.
\end{split}
\end{equation*}
on $M^{n}\times [0,T]$.
\end{lem}

 $\hfill\#$

We remark that the fact that the curvatures of
$h^{a}_{\alpha\beta}$ and $g^{a}_{ij}(\cdot,t)$ are uniformly
bounded (independent of $a$) is essential in our argument. While
the injectivity radius bound $i_a$ may depend on $a$.

For the new family of metrics $g^{a}_{ij}(\cdot,t)$, we have the
following lemma.
\begin{lem}
\begin{eqnarray*}
&\frac{\partial }{\partial
t}g^{a}_{ij}=e^{\varphi^{a}}(-2\overset{a}{R}_{ij}+(\overset{a}{\nabla}^{2}\varphi^{a}+\overset{a}{\nabla}\varphi^{a}\ast
\overset{a}{\nabla}\varphi^{a})\ast
\overset{a}{g}\ast(\overset{a}{g})^{-1}),\\
&\frac{\partial}{\partial t}\overset{a}{\Gamma}^{k}_{ij}=
e^{\varphi^{a}}(\overset{a}{g})^{-1}\ast
\overset{a}{\nabla}\overset{a}{Ric}+e^{\varphi^{a}}(\overset{a}{g})^{-2}\ast\overset{a}{g}\ast
(\overset{a}{Ric}\ast \nabla
\varphi^{a}+\overset{a}{\nabla}^{3}\varphi^{a})\\&\ \
+e^{\varphi^{a}}(\overset{a}{g})^{-3}\ast(\overset{a}{g})^{2}\ast
[(\overset{a}{\nabla}\overset{a}{\varphi})^{3}+\overset{a}{\nabla}^{3}\varphi^{a}],
\end{eqnarray*}
\begin{equation*}\tag{2.20}
\begin{split}
e^{\frac{\varphi^{a}}{2}}|\overset{a}{\nabla}{\varphi}^{a}|_{g^{a}(\cdot,t)}
+e^{\varphi^{a}}|\overset{a}{\nabla}^{2}_{g^{a}(\cdot,t)}\varphi^{a}|_{g^{a}(\cdot,t)}&\leqslant
C(n,k_{0},T),\\
e^{\frac{3}{2}\varphi^{a}}|\overset{a}{\nabla}^{3}_{g^{a}(\cdot,t)}\varphi^{a}|_{g^{a}(\cdot,t)}&\leqslant
C(n,k_{0},T)(1+|\log t|),\\
e^{\frac{k}{2}\varphi^{a}}|\overset{a}{\nabla}^{k}_{g^{a}(\cdot,t)}\varphi^{a}|_{g^{a}(\cdot,t)}&\leqslant
C(n,k_{0},T,k)\frac{1}{t^{\frac{k-3}{2}}},\ \  \mbox{for}\ \
k\geqslant4.
\end{split}
\end{equation*}
\end{lem}
Proof. Note that \begin{equation*}
\begin{split}
\overset{a}{\Gamma}-\Gamma&=g\ast g^{-1}\ast \nabla\varphi^{a}\\
\overset{a}{\nabla^{2}}\varphi^{a}&=\nabla^{2}\varphi^{a}+(\overset{a}{\Gamma}-\Gamma)\ast
\nabla \varphi^{a}\\
\overset{a}{\nabla^{k}}\varphi^{a}&=\sum_{i_1+\cdots+i_p=k}
g^{k-1}\ast (g^{-1})^{k-1}\ast
\nabla^{i_1}\varphi^{a}\ast\cdots\ast\nabla^{i_p}\varphi^{a}
\end{split}
\end{equation*}
where the summation is taken over all indices $i_j>0$. By
combining this with (2.18), we get the desired estimates for
$|\overset{a}{\nabla^{k}}_{g^{a}(\cdot,t)}\varphi^{a}|_{g^{a}(\cdot,t)}$.
One the other hand, since
$$
\overset{a}{R}_{ij}=R_{ij}+(\overset{a}{\nabla^{2}}\varphi^{a}+\nabla\varphi^{a}\ast\nabla\varphi^{a})\ast
g\ast g^{-1},
$$
it follows that
 \begin{equation*}
\begin{split}
\overset{a}{\nabla}_{i}\overset{a}{R}_{jl}&=\nabla_{i}R_{jl}+g^{a}\ast
{g^{a}}^{-1}\ast(\overset{a}{Ric}\ast\nabla\varphi^{a}+\overset{a}{\nabla^{3}}\varphi^{a})\\
&\ \ +{(g^{a})}^{2}\ast {(g^{a})}^{-2}\ast
(\overset{a}{\nabla^{2}}\varphi^{a}\ast\overset{a}{\nabla}\varphi^{a}+(\overset{a}{\nabla}\varphi^{a})^{3}).
\end{split}
\end{equation*}
By combining this with
\begin{equation*}
\begin{split}
\frac{\partial}{\partial
t}\overset{a}{\Gamma}^{k}_{ij}&=\frac{\partial}{\partial
t}{\Gamma}^{k}_{ij}+\frac{\partial}{\partial t}(g^{-1}\ast g\ast
\nabla \varphi^{a})\\
&=-g^{kl}(\nabla_{i}R_{jl}+\nabla_{j}R_{li}-\nabla_{l}R_{ij})+g\ast
g^{-2}\ast Ric \ast \nabla\varphi^{a},
\end{split}
\end{equation*}
we have proved the lemma.

 $\hfill\#$

\subsection{Modified harmonic map flow}
The purpose of this subsection is to solve the equation
$(2.3)_{a}$. More precisely, we will prove the following theorem
\begin{thm}
     There exists $0<T_1<T$, depending only on $k_{0}$, $T$ and $n$ such that for all $a\geqslant1$
     the modified harmonic map flow flow
   coupled with the Ricci flow
    $$ \left\{
     \begin{split}
     \quad \frac{\partial}{\partial t}\overset{a}{F}(x,t)&=\overset{a}{\triangle}_{t} \overset{a}{F}(x,t)  \\
  \overset{a}{F}(\cdot,0)&= identity \  \  \  \  \
  \end{split}
  \right.\eqno(2.3)_{a}
  $$
     has a solution on $M^{n}\times[0,T_0]$ satisfying the
     following estimates
     $$
     \begin{array}{lr}
     |\overset{a}{\nabla}\overset{a}{F}|\leq C(n,k_{0},T),\\
           |{\overset{a}{\nabla}}^{k}\overset{a}{ F}|\leq C(n,k_{0},T,k) t^{-\frac{k-2}{2}},\ \ \ \ \
      \mbox{for all k }\ \geq 2,
     \end{array}\eqno{(2.21)}
     $$
     for some constants $C(n,k_{0},T,k)$ depending only on $n$, $k_{0}$, $T$, and $k$ but independent of $a$.
     \end{thm}

     Note that $\overset{a}{F}$ is viewed as a map from
     $(M^{n},g^{a}_{ij}(x,t))$ and
     $(N^{n},h^{a}_{\alpha\beta}(y))$, all the covariant
     derivatives and the norms in  Theorem 2.6 are computed with
     respect $g^{a}_{ij}(x,t)$ and
     $h^{a}_{\alpha\beta}(y)$. We begin with a easier short time existence of
     $(2.3)_{a}$ where the short time interval may depend on $a$.

\subsubsection{Short time existence of the modified harmonic map flows}
We consider $(2.3)_{a}$ with general initial data.
\begin{thm}
Let $f$ be a smooth map from $M^{n}$ to $N^{n}$ with
$$
E_{0}=\sup_{x\in
M^{n}}|\overset{a}{\nabla}f|_{g^{a}_{ij}(\cdot,0),h^{a}_{\alpha\beta}}(x)+\sup_{x\in
M^{n}}|\overset{a}{\nabla^{2}}f|_{g^{a}_{ij}(\cdot,0),h^{a}_{\alpha\beta}}(x)<\infty.
$$
Then there exists a $\delta_{0}>0$ such that the initial problem
$$ \left\{
     \begin{split}
     \quad \frac{\partial}{\partial t}\overset{a}{F}(x,t)&=\overset{a}{\triangle}_{t} \overset{a}{F}(x,t),  \\
  \overset{a}{F}(x,0)&= f(x), \  \  \  \  \
  \end{split}
  \right.\eqno{(2.3)_{a}}^{\prime}
  $$
  has a smooth solution on $M^{n}\times[0,\delta_{0}]$ satisfying
  the following estimates
$$
\sup_{(x,t)\in
M^{n}\times[0,\delta_0]}|\overset{a}{\nabla}\overset{a}{F}|_{g^{a}_{ij}(\cdot,0),h^{a}_{\alpha\beta}}(x,t)+\sup_{(x,t)\in
M^{n}\times[0,\delta_0]}|\overset{a}{\nabla^{2}}\overset{a}{F}|_{g^{a}_{ij}(\cdot,0),h^{a}_{\alpha\beta}}(x,t)\leqslant
C(n,k_{0},T,a,E_{0}),$$
 $$\sup_{(x,t)\in
M^{n}\times[0,\delta_0]}|\overset{a}{\nabla^{k}}\overset{a}{F}|_{g^{a}_{ij}(\cdot,0),h^{a}_{\alpha\beta}}(x,t)\leqslant
\frac{C(n,k_{0},T,k,a)}{t^{\frac{k-2}{2}}}, \eqno(2.22)$$ for $k
\geq 3.$
\end{thm}

     We will prove the theorem by solving the corresponding
     initial-boundary value problem on a sequence of exhausted bounded
     domains $D_1\subseteq D_2 \subseteq \cdots $ with smooth boundaries and $D_j\supseteq B_{g^a(\cdot,0)}^{a}(P,j+1)$ :
      $$
      \left\{
      \begin{split}
     \quad \frac{\partial}{\partial t}\overset{a}{F^{j}}(x,t)&=\overset{a}{\triangle}_{t} \overset{a}{F^{j}}(x,t),
     \ \ \mbox{ for  }x \in D_{j} \mbox{ and } t>0,\\
  F^{j}(x,0)&= f(x) \  \  \  \  \   \mbox{ for  }x \in D_{j} ,\\
  \overset{a}{F^{j}}(x,t)&= f(x) \ \ \ \ \ \ \ \ \mbox{ for } x\in \partial
  D_{j},
  \end{split}
  \right.\eqno(2.23)
     $$
     and $\overset{a}{F}$ will be obtained as the limit of a convergent subsequence of $\overset{a}{F^{j}}$ as $j\rightarrow
     \infty$. Here $P$ is a fixed point on $M^n$ and
     $B_{g^a(\cdot,0)}^{a}(P,j+1)$ is the geodesic ball centered
     at $P$ of radius $j+1$ with respect to the metric $g^a_{ij}(\cdot,0)$

     The following lemma gives the zero-order estimate of
     $\overset{a}{F^{j}}$.
     \begin{lem}
     There exist positive constants $0<T_2<T$ and $C>0$ such that for any $j$, if
     (2.23)
     has a smooth solution $\overset{a}{F^{j}} $ on
     $\bar{D_{j}}\times[0,T_3]$ with $T_3\leq T_2$, then we have
     $$d_{(N^{n},h^{a})}(f(x),\overset{a}{F^{j}}(x,t))\leq C\sqrt{t}, \eqno(2.24)$$
     for any $(x,t)\in
     D_{j}\times[0,T_3]$.
     \end{lem}
     Proof.  For simplicity, we drop the superscripts $a$ and $j$ of $\overset{a}{F^{j}}$.
     Note that the distance function $d_{(N^{n},h^{a})}(y_{1},y_{2})$ can be regarded as a function on $N^{n}\times
  N^{n}$. Set $\psi(y_{1},y_{2})=\frac{1}{2}d^2_{(N^{n},h^{a})}(y_{1},y_{2})$
  and $\rho(x,t)=\psi(f(x),F(x,t))$. Then $\psi(x,t)$ is smooth
  when $\psi<\frac{1}{2}i_{a}^{2}$.
   Now we compute the equation of $\rho(x,t)$:
$$ (\frac{\partial}{\partial
     t}-\overset{a}{\triangle_{t}})\rho=-d_{h^{a}}(f(x),F(x,t))
     \frac{\partial d_{h^{a}}}{\partial {y_{1}}^{\alpha}}\overset{a}
     {\triangle_{t}}f^{\alpha}-Hess(\psi)(X_{i},X_{j})(g^{a})^{ij}\eqno(2.25)
$$ where the vector fields $X_{i}$, $i=1, 2, \cdots, n$,
in local coordinates $(y_{1}^{\alpha},y_{2}^{\beta})$ on
$N^{n}\times N^{n}$ are defined as follows
$$
X_{i}=\frac{\partial f^{\alpha}}{\partial
x^{i}}\frac{\partial}{\partial y_{1}^{\alpha}}+\frac{\partial
F^{\beta}}{\partial x^{i}}\frac{\partial}{\partial y_{2}^{\beta}}.
$$

To handle the first term on the right hand side of (2.25), we use
$$\overset{a}{\Gamma^{k}_{ij}}(x,t)-\overset{a}{\Gamma^{k}_{ij}}(x,0)=
\Gamma^{k}_{ij}(x,t)-{\Gamma}^{k}_{ij}(x,0)+g(\cdot,t)\ast
g^{-1}(\cdot,t)\ast \nabla \varphi^{a}+g(\cdot,0)\ast
g^{-1}(\cdot,0)\ast \nabla \varphi^{a},$$ to conclude that
$$|\overset{a}{\triangle_{t}}f|_{g^{a}(\cdot,t),h^{a}}\leqslant C(n,k_{0},T)E_{0}.$$
  Recall from Lemma 2.3 that the curvature of the metric $h^a_{\alpha\beta}$ is bounded by $\bar{C_0}$.
  We claim that if $d_{h^{a}}(f(x),F(x,t))\leq \min\{\frac{i_a}{4},\frac{\pi}{4\sqrt{\bar{C_{0}}}}\}$,
  then
  $$
  Hess(\psi)(X_{i},X_{j})(g^{a})^{ij}\geq
  \frac{1}{2}|\overset{a}{\nabla}\overset{a}{F}|_{g^{a},h^{a}}^{2}-C \eqno(2.26)
  $$
  where $C=C(E_0,\bar{C_0})$ depends only on $E_0$ and $\bar{C_0}$.

  Indeed, recall the computation of $Hess(\psi)$ in \cite{ScY1}. For any $(u,v)\in
D=\{(u,v): (u,v)\in N^{n}\times
N^{n},d_{N^{n}}(u,v)<\min\{\frac{i_{a}}{2},\frac{\pi}{2\sqrt{\bar{C_0}}}\}\}$,
let $\gamma_{uv}$ be the minimal geodesic from $u$ to $v$ and
$e_{1}\in T_{u}N^{n}$ be the tangent vector to $\gamma_{uv}$ at
$u$. Then $e_{1}(u,v)$ defines a smooth vector field on $D$.  Let
$\{e_{i}\}$ be an orthonormal basis for $T_{u}N^{n}$ which depends
$u$ smoothly. By parallel translation of $\{e_{i}\}$ along
$\gamma$, we define $\{\bar{e}_{i}\}$ an orthonormal basis
 for $T_{v}N^{n}$. Thus $\{e_{1},\cdots e_{n},\bar{e}_{1},\cdots \bar{e}_{n}\}$
 is a local frame on $D$. Then For any $X=X^{(1)}+X^{(2)}\in T_{(u,v)}D$, where
$$
X^{(1)}=\sum_{i=1}^{n}\xi_{i}e_{i}, \ \  \mbox{and}\ \
X^{(2)}=\sum_{i=1}^{n}\eta_{i}\bar{e}_{i},
$$
  by the formula (16) in
  \cite{ScY1},
\begin{eqnarray*}
Hess(\psi)(X,X)&=&\sum_{i=1}^{n}(\xi_{i}-\eta_i)^{2}
+\int_{0}^{r}t\langle\nabla_{e_{1}}V,\nabla_{e_{1}}V\rangle
+\int_{0}^{r}t\langle\nabla_{\bar{e}_{1}}V,\nabla_{\bar{e}_{1}}V\rangle\\
&&-\int_{0}^{r}t\langle
R(e_{1},V)V,e_{1}\rangle-\int_{0}^{r}t\langle
R(\bar{e}_{1},V)V,\bar{e}_{1}\rangle
\end{eqnarray*}
where $V$ is a Jacobi field on geodesic $\sigma$ (connecting
$(v,v)$ to $(u,v)$) and $\bar{\sigma}$ (connecting $(u,u)$ to
$(u,v)$) with $X$ as the boundary values, where $X$ is extended to
be a local vector field by letting its coefficients with respect
to $\{e_{1},\cdots e_{n},\bar{e}_{1},\cdots \bar{e}_{n}\}$ be
constant(see \cite{ScY1}). By the Jacobi equation, $|V|$,
$|\nabla_{e_{1}}V|$ and  $|\nabla_{\bar{e}_{1}}V|$ are bounded.
Thus
   we have $$ |Hess(\psi)|_{h^{a}}\leq C(i_a,\bar{C_0})$$ under
  the assumption of the claim. So the mixed term
  $\frac{\partial^{2}\psi}{\partial y^{\alpha}_{1}\partial y^{\beta}_{2}}f^{\alpha}_{i}F^{\beta}_{j}(g^{a})^{ij}$
  in $Hess(\psi)(X_i,X_j)(g^{a})^{ij}$ can be bounded by
  $C(E_0,\bar{C_0})E_0|\overset{a}{\nabla}\overset{a}{F}|_{g^{a},h^{a}}$. On the other hand, the Hessian
  comparison theorem for the points which are not in the cut locus gives
  \begin{eqnarray*}
   \frac{\partial \psi}{\partial {y_{2}}^{\alpha}\partial
    {y_{2}}^{\beta}}
     -(\overset{a}{\Gamma^{\gamma}_{\alpha\beta}}\circ \overset{a}{F})\frac{\partial \psi}{\partial
     y_{2}^{\gamma}}\geq \frac{\pi}{4}{h}^{a}_{\alpha\beta},\\
\frac{\partial \psi}{\partial {y_{1}}^{\alpha}\partial
{y_{1}}^{\beta}}
     -(\overset{a}{\Gamma^{\gamma}_{\alpha\beta}}\circ f)\frac{\partial \psi}{\partial
     y_{1}^{\gamma}}\geq \frac{\pi}{4}{h}^{a}_{\alpha\beta}.
  \end{eqnarray*}
Thus the claim follows.

  Let $$T^{\prime}_2=\max \{t\leq T: \sup\limits_{D} {d}_{h^{a}}(f(x),F(x,t))\leq
  \min\{i_a,\frac{\pi}{4\sqrt{\bar{C_0}}}\}\}.$$
   If $\overset{a}{F}(x,t)$ is a smooth solution of (2.23)
   on $\bar{D}\times[0,T_3]$ with $T_3\leqslant T^{\prime}_{2}$, by (2.25) and (2.26), we get
   \begin{equation*} \tag{2.27}     \begin{split}
(\frac{\partial}{\partial
     t}-\overset{a}{\triangle_{t}})\rho\leq -\frac{1}{2}|\overset{a}{\nabla}\overset{a}{F}|_{g^{a},h^{a}}^{2}
     +C\sqrt{\rho}
     +C
     \end{split}
\end{equation*}
on $D\times[0,T_3]$, for some constant $C$ depending on $E_0$,
$i_{a}$ and $\bar{C_{0}}$. Note that the initial and boundary
values of $\rho$ are zero, so by the maximum principle, we get
$$
 d_{h^{a}}(f(x),\overset{a}{F}(x,t))\leq C\sqrt{t}.
$$
This implies $$T^{\prime}_2\geq
\min\{\frac{\min\{i_{a},\frac{\pi}{4\sqrt{\bar{C_{0}}}}\}^{2}}{C^{2}},T_3\}.$$
Hence the lemma holds with
$$T_2=\min\{\frac{\min\{i_{a},\frac{\pi}{4\sqrt{\bar{C_{0}}}}\}^{2}}{C^{2}},T\}.$$
 .

 $\hfill\#$

 After we have the zero order estimate (2.24), we now apply the standard parabolic equation
 theory to get the following short time existence for (2.23).
 \begin{lem}
  There exists a positive constant $T_3 \leq  T_2$ depending
 only on the dimension $n, a, T_2$ and $C$ in Lemma 2.8 such that for each $j$,
 the initial-boundary value problem (2.23) has a smooth
 solution $\overset{a}{F^{j}}$ on $\bar{D_{j}}\times[0,T_3]$.
 \end{lem}

 $\underline{\mbox{\textbf{Proof}}}$. For an arbitrarily fixed point $x_0$ in $\bar{D_j}$,
 choose normal coordinates $\{x^i\}$ and $\{y^\alpha \}$ on $(M^{n},g^{a}(\cdot,0))$ and $(N^{n},h^{a})$
 around $x_0$ and $f(x_0)$ respectively. The equation (2.23) can be written as
   $$ \frac
{\partial {y^{\alpha}}}{\partial t} (x,t) = (g^{a})^{ij}(x,t)
\{\frac{\partial^{2}y^{\alpha}}{\partial
    x^{i}\partial
   x^{j}}-\overset{a}{\Gamma^{k}_{ij}}(x,t)\frac{\partial y^{\alpha}}{\partial
   x^{k}}
   +\overset{a}{\Gamma^{\alpha}_{\beta\gamma}}(y^1(x,t), \cdots, y^n(x,t))\frac{\partial y^{\beta}}{\partial x^{i}}
   \frac{\partial y^{\gamma}}{\partial x^{j}}\}.\eqno(2.28)$$
   Note that $\overset{a}{\Gamma^{\alpha}_{\beta\gamma}}(f(x_0))=0$. By applying (2.24) and a result of Hamilton
   (Corollary (4.12) in \cite{Ha4}), we
   know that the coefficients of the
   quadratic terms on the RHS of (2.28) can be as small as we
   like provided $T_3>0$ sufficiently small (independent of
   $x_0$ and $j$).

   Now for fixed $j$, we consider the corresponding parabolic
   system of the difference of the map $\overset{a}{F^j}$ and $f(x)$.
   Clearly the coefficients of the quadratic terms of the
   gradients are also very small. Thus, whenever (2.23) has a
   solution on a time interval $[0,T'_3]$ with $T'_3 \leq T_3$, we
   can argue exactly as in the proof of Theorem 6.1 in Chapter VII
   of the book \cite{LSU} to bound the norm of $\overset{a}{\nabla}
    \overset{a}{F^{j}}$ over $\bar{D_j}\times[0,T'_3]$ by a constant depending only
   on the $L^{\infty}$ bound of $\overset{a}{F}$ in (2.23), the map $f(x)$,
   the domain $D_j$, and the metrics $g^{a}_{ij}(\cdot,t)$ and $h^{a}_{\alpha
   \beta}$ over the domain $D_{j+1}$. Hence  by the same argument as in the proof of
   Theorem 7.1 in Chapter VII of the book \cite{LSU}, we deduce that the initial-boundary value problem (2.23) has a smooth
 solution $\overset{a}{F^{j}}$ on $\bar{D_{j}}\times[0,T_3]$.

 $\hfill\#$

 Unfortunately, the gradient estimates of $\overset{a}{F^{j}}$ in
 the proof of the above lemma depend also on the domain $D_j$. In order to get a convergent subsequence of
 $\overset{a}{F^{j}}$, we have to estimate the covariant
 derivatives of $\overset{a}{F^{j}}$ uniformly in each compact subsets. Before we proceed,
 we need some preliminary estimates.
 \begin{lem}
 The covariant derivatives of $\overset{a}{F^{j}}$ satisfy the
 following equations
\begin{equation*} \tag{2.29}     \begin{split}
     \frac{\partial}{\partial t}\overset{a}{\nabla}\overset{a}{F^{j}}
     &=\overset{a}{\triangle}_{t} \overset{a}{\nabla}\overset{a}{F^{j}}+
      \overset{a}{Ric}(M^{n})\ast\overset{a}{\nabla}\overset{a}{F^{j}}
      +\overset{a}{R_{N}}\ast(\overset{a}{\nabla}\overset{a}{F^{j}})^{3}, \\
    \frac{\partial}{\partial t}\overset{a}{\nabla^{k}}\overset{a}{F^{j}}
     &=\overset{a}{\triangle}_{t} \overset{a}{\nabla^{k}}\overset{a}{F^{j}}+
      \sum_{l=0}^{k-1}\overset{a}{\nabla^{l}}[(\overset{a}{R_M}+\overset{a}{R_{N}}\ast
      (\overset{a}{\nabla}\overset{a}{F^{j}})^{2}+e^{\varphi^{a}}\overset{a}{R_M}
      +\overset{a}{\nabla^{2}}e^{\varphi^{a}})\ast\overset{a}{\nabla^{k-l}}\overset{a}{F^{j}}],
  \end{split}
            \end{equation*}
            where $\overset{a}{\nabla^{l}}(A\ast B)$ represents
            the linear combinations of $\overset{a}{\nabla^{l}}A\ast
            B$,$\overset{a}{\nabla^{l-1}}A\ast
            \overset{a}{\nabla}B$, $\cdots$, $A\ast
            \overset{a}{\nabla^{l}}B$, and
            $\overset{a}{\nabla^{2}}e^{\varphi^{a}}=e^{\varphi^{a}}(\overset{a}{\nabla^{2}}
            \varphi^{a}+\overset{a}{\nabla}\varphi^{a}\ast \overset{a}{\nabla}\varphi^{a}).$
 \end{lem}
 $\underline{\mbox{\textbf{Proof}}}$
 For $k=1$, by direct computation and Ricci formula, we have
$$
\frac{\partial}{\partial
t}\overset{a}{\nabla_{i}}\overset{a}{F^{\alpha}}=\overset{a}{\triangle_{t}}\overset{a}{\nabla_{i}}\overset{a}{F^{\alpha}}
-\overset{a}{R^{l}_{i}}\overset{a}{\nabla_{l}}\overset{a}{F^{\alpha}}+\overset{a}{R}^{\alpha}_{\beta\delta\gamma}
\overset{a}{\nabla_{i}}\overset{a}{F^{\beta}}\overset{a}{\nabla_{k}}\overset{a}{F^{\delta}}\overset{a}{\nabla_{l}}\overset{a}{F^{\gamma}}(
g^{a})^{kl}.
$$
For $k\geqslant2$, by Ricci formula, it follows
$$
\overset{a}{\nabla}\overset{a}{\triangle}\overset{a}{\nabla^{k-1}}\overset{a}{F^{j}}=\overset{a}{\triangle}
\overset{a}{\nabla^{k}}\overset{a}{F^{j}}+\overset{a}{\nabla}
[(\overset{a}{R_M}+\overset{a}{R_N}\ast(\overset{a}{\nabla}\overset{a}{F^{j}})^{2})\ast
\overset{a}{\nabla^{k-1}}\overset{a}{F^{j}}].
$$
Recall from (2.20) that
$$
\frac{\partial}{\partial
t}\overset{a}{\Gamma^{i}_{jk}}=\overset{a}{\nabla}(e^{\varphi^{a}}\overset{a}{R_M}+\overset{a}{\nabla^{2}}e^{\varphi^{a}}).
$$
Then we have
\begin{equation*}      \begin{split}
    \frac{\partial}{\partial t}\overset{a}{\nabla^{k}}\overset{a}{F^{j}}
     -\overset{a}{\triangle}_{t}
     \overset{a}{\nabla^{k}}\overset{a}{F^{j}}=&\overset{a}{\nabla}[(\frac
     {\partial}{\partial t}-\overset{a}{\triangle_{t}})\overset{a}{\nabla^{k-1}}
     \overset{a}{F^{j}}]+\overset{a}{\nabla}(e^{\varphi^{a}}\overset{a}{R_M}
     +\overset{a}{\nabla^{2}}e^{\varphi^{a}})
     \ast\overset{a}{\nabla^{k-1}}\overset{a}{F^{j}}\\&\ \
     +\overset{a}{R_N}\ast\overset{a}{\nabla}\overset{a}
     {F^{j}}\ast\overset{a}{\nabla^{2}}\overset{a}{F^{j}}\ast
     \overset{a}{\nabla^{k-1}}\overset{a}{F^{j}}+\overset{a}
     {\nabla}[(\overset{a}{R_M}+\overset{a}{R_N}\ast(\overset{a}
     {\nabla}\overset{a}{F^{j}})^{2})\ast \overset{a}{\nabla^{k-1}}\overset{a}{F^{j}}]\\
     =& \overset{a}{\nabla}[(\frac{\partial}{\partial t}
     -\overset{a}{\triangle})\overset{a}{\nabla^{k-1}}
     \overset{a}{F^{j}}]+\overset{a}{\nabla}\{
     (\overset{a}{R_M}+\overset{a}{R_{N}}\ast
      (\overset{a}{\nabla}\overset{a}{F^{j}})^{2}+(e^{\varphi^{a}}\overset{a}{R_M}
      +\overset{a}{\nabla^{2}}e^{\varphi^{a}})\ast\overset{a}{\nabla^{k-1}}\overset{a}{F^{j}}\}\\
     =&
      \sum_{l=0}^{k-1}\overset{a}{\nabla^{l}}[(\overset{a}{R_M}+\overset{a}{R_{N}}\ast
      (\overset{a}{\nabla}\overset{a}{F^{j}})^{2}+e^{\varphi^{a}}\overset{a}{R_M}
      +\overset{a}{\nabla^{2}}e^{\varphi^{a}})\ast\overset{a}{\nabla^{k-l}}\overset{a}{F^{j}}].
  \end{split}
            \end{equation*}
            This proves the lemma.

            $\hfill\#$

For each $k>0$, let $\xi_k$ be a smooth non-increasing function
from $(-\infty,+\infty)$ to $[0,1]$ so that $\xi_k(s)=1$ for
   $s\in (-\infty,\frac{1}{2}+\frac{1}{2^{k+1}}]$, and $\xi_k(s)=0$ for $s\in[\frac{1}{2}+\frac{1}{2^k})$; moreover for
   any $\epsilon>0$ there exists a universal $C_{k,\epsilon}>0$ such
   that
   $$
|\xi^{\prime}_k(s)|+|\xi^{\prime\prime}_k(s)|\leq
C_{k,\epsilon}{\xi_k(s)}^{1-\epsilon}.
$$
  \begin{lem} There exists a positive constant $T_4$, $0<T_4\leq T_3$ independent of $j$ such that for any geodesic ball
  $B_{g^{a}(\cdot,0)}(x_0,\delta)\subset D_j$, there is a constant
  $C=C(a,\delta,E_0,\bar{C_0},\bar{k_0})$ such that the smooth
  solution of (2.23) satisfies
$$   |\overset{a}{\nabla}\overset{a}{
     F^{j}}|_{g^{a}(\cdot,t),h^{a}}\leqslant C
     $$
    on $B_{g^{a}(\cdot,0)}(x_0,\frac{3\delta}{4})\times [0,T_4]$.
\end{lem}
$\underline{\mbox{\textbf{Proof.}}}$ We compute the equation of
$|\overset{a}{\nabla}\overset{a}{
     F^{j}}|^{2}_{g^{a}(\cdot,t),h^{a}}.$ For simplicity, we drop the superscript $j$. By (2.20), we have
\begin{equation*} \tag{2.30}     \begin{split}
     (\frac{\partial}{\partial
     t}-\overset{a}{\triangle}_{t})|\overset{a}{\nabla}
     \overset{a}{F}|^{2}_{g^{a}(\cdot,t),h^{a}}
     &=\langle \overset{a}{Ric}(M^{n})\ast\overset{a}
     {\nabla}\overset{a}{F}
      +\overset{a}{R_{N}}\ast(\overset{a}{\nabla}\overset{a}{F})^{3},
      \overset{a}{\nabla}\overset{a}{F}\rangle_{g^{a},h^{a}}
      -2|\overset{a}{\nabla^{2}}
      \overset{a}{F}|^{2}_{g^{a}(\cdot,t),h^{a}}\\
      &\ \ +e^{\varphi^{a}}
      (\overset{a}{Ric(M^{n})}+\overset{a}{\nabla^{2}}
      \varphi^{a}+
      \overset{a}{\nabla}\varphi^{a}\ast\overset{a}
      {\nabla}\varphi^{a})
      \ast\overset{a}{\nabla}
      \overset{a}{F}\ast\overset{a}{\nabla}\overset{a}{F}\\&\leqslant
      -2|\overset{a}{\nabla^{2}}\overset{a}{F}|^{2}_{g^{a}(\cdot,t),h^{a}}+C(n,k_0,T)|\overset{a}{\nabla}
      \overset{a}{F}|^{2}_{g^{a}(\cdot,t),h^{a}}+C(n)\bar{C_0}
      |\overset{a}{\nabla}\overset{a}{F}|^{4}_{g^{a}(\cdot,t),h^{a}}.
        \end{split}
       \end{equation*}
       Setting
       $$
       \rho_{A}(x,t)=(d_{h^{a}}^{2}(f(x),F(x,t))+A)|\overset{a}{\nabla}\overset{a}{
     F}|^{2}_{g^{a}(\cdot,t),h^{a}}
       $$
       where $A$ is determined later, and combining with (2.27) and (2.24),
       we have
\begin{eqnarray*}
\frac{\partial}{\partial
t}\rho_{A}&\leqslant&\overset{a}{\triangle}\rho_A-2|\overset{a}{\nabla^{2}}
\overset{a}{F}|^{2}_{g^{a},h^{a}}(d_{h^{a}}^{2}(f(x),\overset{a}{F}(x,t))+A)-
|\overset{a}{\nabla}\overset{a}{ F}|^{4}_{g^{a},h^{a}}\\&
&+C(n)\bar{C_0}(d^{2}_{h^{a}}(f(x),\overset{a}{F}(x,t))+A)|\overset{a}{\nabla}
\overset{a}{F}|^{4}_{g^{a},h^{a}}+C|\overset{a}{\nabla}\overset{a}{F}|^{2}_{g^{a},h^{a}}+C(n,k_0,T)\rho_A
\\& &+2|\nabla
d^{2}_{h^{a}}(f(x),\overset{a}{F}(x,t))|_{g^{a}} |\nabla
|\overset{a}{\nabla}\overset{a}{
     F}|^{2}_{g^{a}(\cdot,t),h^{a}}|_{g^{a}}.
 \end{eqnarray*}
Since
\begin{eqnarray*}
|\nabla
d^{2}_{h^{a}}(f(x),\overset{a}{F}(x,t))|_{g^{a}}&\leqslant& 2
d_{h^{a}}(f(x),\overset{a}{F}(x,t))(|\overset{a}{\nabla}\overset{a}{F}|_{g^{a},h^{a}}+
|\overset{a}{\nabla}f|_{g^{a},h^{a}})\\
&\leqslant&
C\sqrt{t}+C\sqrt{t}|\overset{a}{\nabla}\overset{a}{F}|_{g^{a},h^{a}},\\
 |\nabla
|\overset{a}{\nabla}\overset{a}{
     F}|^{2}_{g^{a}(\cdot,t),h^{a}}|_{g^{a}}&\leqslant& 2|\overset{a}{\nabla^{2}}\overset{a}{F}|
     _{g^{a},h^{a}}|\overset{a}{\nabla}\overset{a}{F}|_{g^{a},h^{a}},
 \end{eqnarray*}
 by choosing $T_4=\min\{T_3,\frac{1}{4C(n)\bar{C_0}C^{2}}\}$,
 $A=\frac{1}{4C(n)\bar{C_0}}$, and applying Cauchy-Schwartz
 inequality, we have
 $$
(\frac{\partial}{\partial
t}-\overset{a}{\triangle})\rho_{A}\leqslant-(C(n)\bar{C_0})\rho_A^{2}+C.
 $$
Here and in the following we denote by $C$ various constants
depending only on $n$, $k_0$, $T$, $E_0$ and $a$.

We compute the equation of
$u=\xi_{1}(\frac{d_{g^{a}(\cdot,0)}(x_0,\cdot)}{\delta})\rho_{A}$
at the  smooth points of  $d_{g^{a}(\cdot,0)}(x_0,\cdot)$,
$$
(\frac{\partial}{\partial t}-\overset{a}{\triangle})u\leqslant
C\xi_{1}-(C(n)\bar{C_0})\rho_{A}^{2}\xi_{1}-2(g^{a})^{ij}\nabla_{i}\xi_{1}\nabla_{j}\rho_{A}
+(-\xi_{1}^{\prime}\frac{\overset{a}{\triangle}d_{g^{a}(\cdot,0)}(x_0,\cdot)}
{\delta}+e^{nk_0T}\frac{|\xi_{1}^{\prime\prime}|}{\delta^{2}})\rho_{A}.
$$
By the Hessian comparison theorem and the fact that
$-\xi^{\prime}_{1}\geqslant0$, we have
\begin{eqnarray*}
\overset{a}{\nabla_{i}}\overset{a}{\nabla_{j}}d_{g^{a}(\cdot,0)}&\leqslant&
\overset{a}{\nabla_{i}^{0}}\overset{a}{\nabla_{j}^{0}}d_{g^{a}(\cdot,0)}
+(\overset{a}{\Gamma}(\cdot,0)-\overset{a}{\Gamma}(\cdot,t))\ast
\nabla d_{g^{a}(\cdot,0)}\\
&
\leqslant&(\frac{1+\bar{k_0}d_{g^{a}(\cdot,0)}}{d_{g^{a}(\cdot,0)}}+C)g^{a}_{ij}(\cdot,0),\\
-\xi_{1}^{\prime}\overset{a}{\triangle}d_{g^{a}(\cdot,0)}&\leqslant&\frac{C|\xi_{1}^{\prime}|}{\delta}.
\end{eqnarray*}
These two inequalities hold on the whole manifold in the sense of
support functions. Thus for any $x_{1}\in M^{n}$, there is a
function $h_{x_{1}}$ which is smooth on a neighborhood of $x_{1}$
with $h_{x_{1}}(\cdot)\geqslant d_{g^{a}(\cdot,0)}(x_0,\cdot)$,
$h_{x_{1}}(x_{1})= d_{g^{a}(\cdot,0)}(x_0,x_{1})$ and

$$
-\xi_{1}^{\prime}\overset{a}{\triangle}h_{x_{1}}\mid_{x_{1}}\leqslant2\frac{C|\xi_{1}^{\prime}|}{\delta}.
$$
Indeed, $h_{x_{1}}$ can be chosen to having the form
$d_{g^{a}(\cdot,0)}(q,\cdot)+d_{g^{a}(\cdot,0)}(q,x_{0})$ for some
$q$, so we may require
$|\overset{a}{\nabla}h_{x_{1}}|_{g^{a}(\cdot,0)}\leqslant1$.
 Let
$(x_{1},t_0)$ be the maximum point of $u$ over
$M^{n}\times[0,T_{4}]$. If $t_{0}=0$, then
$\xi_{1}\rho_{A}\leqslant E_{0}$. Assume $t_{0}>0$. At the point
$(x_{1},t_{0})$, we have
 $\frac{\partial}{\partial
 t}(\xi_{1}\rho_{A})(x_{1},t_{0})\geqslant0$. If $x_{1}$ does not lie on the cut locus of $x_{0}$, then
\begin{eqnarray*}
0&\leqslant&-C(n)\bar{C_0}\rho_{A}^{2}\xi_{1}+\frac{1}{\delta^{2}}
(e^{nk_0T}\frac{|\xi_{1}^{\prime}|^{2}}{\xi_{1}}+2C(|\xi_{1}^{\prime}|+|\xi_{1}^{\prime\prime}|))\rho_{A}+C\xi_{1}\\
&\leqslant
&-C(n)\bar{C_0}\rho_{A}^{2}\xi_{1}+\frac{C}{\delta^{2}}\sqrt{\xi_{1}}\rho_{A}+C\xi_{1}\\
&\leqslant&
-C(n)\bar{C_0}\rho_{A}^{2}\xi_{1}+\frac{C}{\delta^{4}}\\
&\leqslant&
-C(n)\bar{C_0}(\rho_{A}\xi_{1})^{2}+\frac{C}{\delta^{4}}.
\end{eqnarray*}
  We get
$$
\xi_{1}\rho_{A}\leqslant\max\{E_{0},\sqrt{\frac{C}{C(n)\bar{C_{0}}\delta^{4}}}\}
$$
for all $(x,t)\in B_{g^{a}(\cdot,0)}(x_0,\delta)\times[0,T_4]$. If
$x_{1}$ lies on the cut locus of $x_{0}$, then by applying the
standard support function technique (see for example \cite{ScY}),
the above maximum principle argument still works. So by the
definition of $\xi_{1}$ and $\rho_{A}$, we have

$$   |\overset{a}{\nabla}\overset{a}{
     F^{j}}|_{g^{a}(\cdot,t),h^{a}}\leqslant \frac{C}{\delta}
     $$
    on $B_{g^{a}(\cdot,0)}(x_0,\frac{3\delta}{4})\times [0,T_4]$.
    The proof of the lemma is completed.

    $\hfill\#$

    The next lemma estimates the higher derivatives in terms of
    the bound of  $|\overset{a}{\nabla}\overset{a}{
     F^{j}}|_{g^{a}(\cdot,t),h^{a}}$.

   \begin{lem} Let $\overset{a}{F}$ be a smooth solution of equation
   $$(\frac{\partial}{\partial
   t}-\overset{a}{\triangle})\overset{a}{F}=0
   $$
   on $B_{g^{a}(\cdot,0)}(x_{0},\delta)\times[0,\bar{T}]$, with $\bar{T}\leqslant T$. Suppose
   \begin{equation*}\tag{2.31}
   \begin{split}
\sup_{(x,t)\in
B_{g^{a}(\cdot,0)}(x_0,\frac{3\delta}{4})\times[0,\bar{T}]}|\overset{a}
{\nabla}\overset{a}{F}|_{g^{a}_{ij}(\cdot,0),h^{a}_{\alpha\beta}}(x,t)&\leqslant
E_{1},\\  \mbox{and  } \sup_{x\in
B_{g^{a}(\cdot,0)}(x_0,\frac{3\delta}{4})}
|\overset{a}{\nabla^{2}}\overset{a}{F}|_{g^{a}_{ij}(\cdot,0),h^{a}_{\alpha\beta}}(x,0)&\leqslant
E_{1}.
 \end{split}
 \end{equation*}
 Then for any $k\geqslant2$, there exists a positive
constant $C=C(k,E_{1},\delta,k_{0},T)>0$ such that
$$ |\overset{a}{\nabla^{k}} \overset{a}{F}|_{g^{a}(\cdot,t),h^{a}}\leq C t^{-\frac{k-2}{2}}
         \eqno(2.32)  $$
    on $B_{g^a(\cdot,0)}(x_{0},\frac{\delta}{2})\times [0,\bar{T}]$.
\end{lem}
Proof. The proof is using the Bernstein trick. We assume
$\delta<1$ without loss of generality. For $k=2$, from (2.15),
(2.19), (2.20) and (2.29), we have
\begin{equation*} \tag{2.33}     \begin{split}
     (\frac{\partial}{\partial
     t}-\overset{a}{\triangle}_{t})|\overset{a}{\nabla^{2}}
     \overset{a}{F}|^{2}_{g^{a}(\cdot,t),h^{a}}
     &=\langle \sum_{l=0}^{1}\overset{a}{\nabla^{l}}[(\overset{a}{R_M}+\overset{a}{R_{N}}\ast
      (\overset{a}{\nabla}\overset{a}{F^{j}})^{2}+e^{\varphi^{a}}\overset{a}{R_M}
      +\overset{a}{\nabla^{2}}e^{\varphi^{a}})\ast\overset{a}{\nabla^{2-l}}\overset{a}{F}],\overset{a}{\nabla^{2}}
     \overset{a}{F}\rangle_{g^{a},h^{a}}\\&\ \ \ \
      -2|\overset{a}{\nabla^{3}}
      \overset{a}{F}|^{2}_{g^{a}(\cdot,t),h^{a}} +e^{\varphi^{a}}
      (\overset{a}{Ric(M^{n})}+\overset{a}{\nabla^{2}}
      \varphi^{a}+
      \overset{a}{\nabla}\varphi^{a}\ast\overset{a}
      {\nabla}\varphi^{a})
      \ast(\overset{a}{\nabla^{2}}
      \overset{a}{F})^{2}\\&\leqslant
      -2|\overset{a}{\nabla^{3}}\overset{a}{F}|^{2}_{g^{a}(\cdot,t),h^{a}}+C|\overset{a}{\nabla^{2}}
      \overset{a}{F}|^{2}_{g^{a}(\cdot,t),h^{a}}+\frac{C}{\sqrt{t}}
      |\overset{a}{\nabla^{2}}\overset{a}{F}|_{g^{a}(\cdot,t),h^{a}}.
        \end{split}
       \end{equation*}
In this lemma, we use $C$ to denote various constants depending
only on $E_1$, $k_0$, $T$, $k$ and $\delta$.  Note that by (2.30)
and (2.33), we have
\begin{equation*}    \begin{split}
     (\frac{\partial}{\partial
     t}-\overset{a}{\triangle}_{t})|\overset{a}{\nabla}
     \overset{a}{F}|^{2}_{g^{a}(\cdot,t),h^{a}}
     &\leqslant
      -2|\overset{a}{\nabla^{2}}\overset{a}{F}|^{2}_{g^{a}(\cdot,t),h^{a}}+C,\\
      (\frac{\partial}{\partial
     t}-\overset{a}{\triangle}_{t})|\overset{a}{\nabla^{2}}
     \overset{a}{F}|_{g^{a}(\cdot,t),h^{a}}
     &\leqslant
      C|\overset{a}{\nabla^{2}}
      \overset{a}{F}|_{g^{a}(\cdot,t),h^{a}}+\frac{C}{\sqrt{t}}.
        \end{split}
       \end{equation*}
So by setting $$v=|\overset{a}{\nabla^{2}}
     \overset{a}{F}|_{g^{a}(\cdot,t),h^{a}}-2C\sqrt{t}+2C\sqrt{T}+|\overset{a}{\nabla}
     \overset{a}{F}|^{2}_{g^{a}(\cdot,t),h^{a}},$$
      we have
     \begin{equation*}     \begin{split}
(\frac{\partial}{\partial
     t}-\overset{a}{\triangle}_{t})v
     &\leqslant
      -2|\overset{a}{\nabla^{2}}\overset{a}{F}|^{2}_{g^{a}(\cdot,t),h^{a}}
      +C|\overset{a}{\nabla^{2}}\overset{a}{F}|_{g^{a}(\cdot,t),h^{a}}+C\\
      &\leqslant -v^{2}+C.
\end{split}
       \end{equation*}
       Since at $t=0$, $$v\leqslant2C\sqrt{T}+E_{1}+E_{1}^{2}$$ on
       $B_{g^{a}(\cdot,0)}(x_{0},\frac{3\delta}{4})$, we
       apply the maximum principle as in Lemma 2.11 to get
       $$
       \xi_{2}(\frac{d_{g^{a}(\cdot,0)}(x_0,\cdot)}{\delta})v\leqslant C
       $$
       on
       $B_{g^{a}(\cdot,0)}(x_{0},\frac{3\delta}{4})\times[0,\bar{T}]$.
       This implies
$$ |\overset{a}{\nabla^{2}} \overset{a}{F}|_{g^{a}(\cdot,t),h^{a}}\leq C
$$
       on
       $B_{g^{a}(\cdot,0)}(x_{0},(\frac{1}{2}+\frac{1}{2^{3}})\delta)\times[0,\bar{T}]$.

Now we estimate the third-order derivatives. From Shi's gradient
estimate \cite{Sh1}, the estimate of $ |\overset{a}{\nabla^{2}}
\overset{a}{F}|_{g^{a}(\cdot,t),h^{a}}\leq C $ and (2.15), (2.19),
(2.20) and (2.29), we have:
\begin{equation*} \tag{2.34}     \begin{split}
     (\frac{\partial}{\partial
     t}-\overset{a}{\triangle}_{t})|\overset{a}{\nabla^{3}}
     \overset{a}{F}|^{2}_{g^{a}(\cdot,t),h^{a}}
     &=\langle \sum_{l=0}^{2}\overset{a}{\nabla^{l}}[(\overset{a}{R_M}+\overset{a}{R_{N}}\ast
      (\overset{a}{\nabla}\overset{a}{F})^{2}+e^{\varphi^{a}}\overset{a}{R_M}
      +\overset{a}{\nabla^{2}}e^{\varphi^{a}})\ast\overset{a}{\nabla^{3-l}}\overset{a}{F}],\overset{a}{\nabla^{3}}
     \overset{a}{F}\rangle_{g^{a},h^{a}}\\&\ \ \ \
      -2|\overset{a}{\nabla^{4}}
      \overset{a}{F}|^{2}_{g^{a}(\cdot,t),h^{a}} +e^{\varphi^{a}}
      (\overset{a}{Ric(M^{n})}+\overset{a}{\nabla^{2}}
      \varphi^{a}+
      \overset{a}{\nabla}\varphi^{a}\ast\overset{a}
      {\nabla}\varphi^{a})
      \ast(\overset{a}{\nabla^{3}}
      \overset{a}{F})^{2}\\&\leqslant
      -2|\overset{a}{\nabla^{4}}\overset{a}{F}|^{2}_{g^{a}(\cdot,t),h^{a}}+C|\overset{a}{\nabla^{3}}
      \overset{a}{F}|^{2}_{g^{a}(\cdot,t),h^{a}}+\frac{C}{t}
      |\overset{a}{\nabla^{3}}\overset{a}{F}|_{g^{a}(\cdot,t),h^{a}}.
        \end{split}
       \end{equation*}
 on
 $B_{g^{a}(\cdot,0)}(x_{0},(\frac{1}{2}+\frac{1}{8})\delta)\times[0,\bar{T}]$.
 Here we used the estimates $|\overset{a}{\nabla^{4}}e^{\varphi^{a}}|_{g^{a}}\leqslant \frac{C}{\sqrt{t}}$, $|\overset{a}{\nabla^{3}}e^{\varphi^{a}}|_{g^{a}}\leqslant C(1+|\log t|)$,
  and $|\overset{a}{\nabla^{2}}\overset{a}{R_m}|\leqslant\frac{C}{t}$.

  By (2.33), it follows
$$(\frac{\partial}{\partial
     t}-\overset{a}{\triangle}_{t})|\overset{a}{\nabla^{2}}
     \overset{a}{F}|^{2}_{g^{a}(\cdot,t),h^{a}}
     \leqslant
      -2|\overset{a}{\nabla^{3}}
      \overset{a}{F}|^{2}_{g^{a}(\cdot,t),h^{a}}+\frac{C}{\sqrt{t}} \eqno(2.35)$$
      on
      $B_{g^{a}(\cdot,0)}(x_{0},(\frac{1}{2}+\frac{1}{8})\delta)\times[0,\bar{T}]$.
Let $v=(|\overset{a}{\nabla^{2}}
      \overset{a}{F}|^{2}_{g^{a}(\cdot,t),h^{a}}+A)|\overset{a}{\nabla^{3}}
      \overset{a}{F}|^{2}_{g^{a}(\cdot,t),h^{a}}$, where $A=100 \sup\limits_{B_{g^{a}(\cdot,0)}(x_0,(\frac{1}{2}+\frac{1}{2^{3}})\delta)\times[0,\bar{T}]}|\overset{a}
{\nabla^{2}}\overset{a}{F}|_{g^{a}_{ij}(\cdot,t),h^{a}_{\alpha\beta}}(x,t)+C.$
By a direct computation, it follows
\begin{eqnarray*}
(\frac{\partial}{\partial
t}-\overset{a}{\triangle})v&\leq&|\overset{a}{\nabla^{3}}
      \overset{a}{F}|^{2}_{g^{a}(\cdot,t),h^{a}}(-2|\overset{a}{\nabla^{3}}
      \overset{a}{F}|^{2}_{g^{a}(\cdot,t),h^{a}}+\frac{C}{\sqrt{t}})+(|\overset{a}{\nabla^{3}}
      \overset{a}{F}|^{2}_{g^{a}(\cdot,t),h^{a}}+A)\\& &\times( -2|\overset{a}{\nabla^{4}}\overset{a}{F}|^{2}_{g^{a}(\cdot,t),h^{a}}+C|\overset{a}{\nabla^{3}}
      \overset{a}{F}|^{2}_{g^{a}(\cdot,t),h^{a}}+\frac{C}{t}
      |\overset{a}{\nabla^{3}}\overset{a}{F}|_{g^{a}(\cdot,t),h^{a}})\\& &
      +
      8|\overset{a}{\nabla^{2}}\overset{a}{F}|_{g^{a}(\cdot,t),h^{a}}
      |\overset{a}{\nabla^{3}}\overset{a}{F}|^{2}_{g^{a}(\cdot,t),h^{a}}
      |\overset{a}{\nabla^{4}}\overset{a}{F}|_{g^{a}(\cdot,t),h^{a}}.
\end{eqnarray*}
Since
 $$8|\overset{a}{\nabla^{2}}\overset{a}{F}|_{g^{a}(\cdot,t),h^{a}}
      |\overset{a}{\nabla^{3}}\overset{a}{F}|^{2}_{g^{a}(\cdot,t),h^{a}}
      |\overset{a}{\nabla^{4}}\overset{a}{F}|_{g^{a}(\cdot,t),h^{a}}\leqslant
      -|\overset{a}{\nabla^{3}}\overset{a}{F}|^{4}_{g^{a}(\cdot,t),h^{a}}+16
      |\overset{a}{\nabla^{4}}\overset{a}{F}|^{2}_{g^{a}(\cdot,t),h^{a}}
      |\overset{a}{\nabla^{2}}\overset{a}{F}|^{2}_{g^{a}(\cdot,t),h^{a}},$$
      we deduce
$$(\frac{\partial}{\partial
t}-\overset{a}{\triangle})v\leqslant-|\overset{a}{\nabla^{3}}\overset{a}{F}|^{4}_{g^{a}(\cdot,t),h^{a}}+
\frac{C}{t}|\overset{a}{\nabla^{3}}\overset{a}{F}|_{g^{a}(\cdot,t),h^{a}}
+\frac{C}{\sqrt{t}}|\overset{a}{\nabla^{3}}\overset{a}{F}|^{2}_{g^{a}(\cdot,t),h^{a}}
$$
and
\begin{equation*}
\begin{split}
(\frac{\partial}{\partial
t}-\overset{a}{\triangle})(tv)&\leqslant
v-t|\overset{a}{\nabla^{3}}\overset{a}{F}|^{4}_{g^{a}(\cdot,t),h^{a}}+
C|\overset{a}{\nabla^{3}}\overset{a}{F}|_{g^{a}(\cdot,t),h^{a}}
+C\sqrt{t}|\overset{a}{\nabla^{3}}\overset{a}{F}|^{2}_{g^{a}(\cdot,t),h^{a}}\\
&\leqslant-\frac{1}{t}\{t^{2}
|\overset{a}{\nabla^{3}}\overset{a}{F}|^{4}_{g^{a}(\cdot,t),h^{a}}-C\sqrt{t}(\sqrt{t}
|\overset{a}{\nabla^{3}}\overset{a}{F}|_{g^{a}(\cdot,t),h^{a}})-tv-C\sqrt{t}(t
|\overset{a}{\nabla^{3}}\overset{a}{F}|^{2}_{g^{a}(\cdot,t),h^{a}})\}\\
&\leqslant-\frac{1}{t}\{\frac{(tv)^{2}}{10^{5}C^{2}}-C\}.
\end{split}
\end{equation*}
So at the maximum point of
$\xi_{3}(\frac{d_{g^{a}(\cdot,0)}(x_{0},\cdot)}{\delta})(tv)$,
applying the maximum principle as in Lemma 2.11, we have
\begin{equation*}
\begin{split}0&\leqslant-\frac{1}{t}\{\frac{\xi_{3}(tv)^{2}}{10^{5}C^{2}}-C\xi_{3}\}
+C(\frac{|\xi_{3}^{\prime}|^{2}}{\xi_{3}}+|\xi_{3}^{\prime\prime}|)(tv)\\
&\leqslant-\frac{1}{t}\{\frac{\xi_{3}(tv)^{2}}{10^{5}C^{2}}-C\xi_{3}
-Ct\sqrt{\xi_{3}}(tv)\}\\
&\leqslant-\frac{1}{t}\{\frac{\xi_{3}(tv)^{2}}{10^{6}C^{2}}-C^{4}\},\\
\end{split}
\end{equation*}
which gives $$ \xi_{3}(tv)\leqslant\sqrt{10^{6}C^{6}}.
$$ Thus by the definition of $v$ and $\xi_{3}$, we get

$$ |\overset{a}{\nabla^{3}} \overset{a}{F}|_{g^{a}(\cdot,t),h^{a}}\leq C t^{-\frac{1}{2}}
         $$
    on $B_{0}(x_{0},(\frac{1}{2}+\frac{1}{2^{4}})\delta)\times [0,\bar{T}]$.

 Now we
 estimate the higher derivatives by induction. Suppose we have
proved that

 $$ |\overset{a}{\nabla^{l}} \overset{a}{F}|_{g^{a}(\cdot,t),h^{a}}\leq C t^{-\frac{l-2}{2}}
        \ \ \ \ \ \  \mbox{ for }
l=3, \cdots, k-1   $$
    on $B_{0}(x_{0},(\frac{1}{2}+\frac{1}{2^{k}})\delta)\times [0,\bar{T}]$. By (2.29), we have
    \begin{equation*}      \begin{split}
     (\frac{\partial}{\partial
     t}-\overset{a}{\triangle}_{t})|\overset{a}{\nabla^{k}}
     \overset{a}{F}|^{2}_{g^{a}(\cdot,t),h^{a}}
     &=\langle \sum_{l=0}^{k-1}\overset{a}{\nabla^{l}}[(\overset{a}{R_M}+\overset{a}{R_{N}}\ast
      (\overset{a}{\nabla}\overset{a}{F})^{2}+e^{\varphi^{a}}\overset{a}{R_M}
      +\overset{a}{\nabla^{2}}e^{\varphi^{a}})\ast\overset{a}{\nabla^{k-l}}\overset{a}{F}],\overset{a}{\nabla^{3}}
     \overset{a}{F}\rangle_{g^{a},h^{a}}\\&\ \ \ \
      -2|\overset{a}{\nabla^{k+1}}
      \overset{a}{F}|^{2}_{g^{a}(\cdot,t),h^{a}} +e^{\varphi^{a}}
      (\overset{a}{Ric(M^{n})}+\overset{a}{\nabla^{2}}
      \varphi^{a}+
      \overset{a}{\nabla}\varphi^{a}\ast\overset{a}
      {\nabla}\varphi^{a})
      \ast(\overset{a}{\nabla^{k}}
      \overset{a}{F})^{2}\\&\leqslant
      -2|\overset{a}{\nabla^{k+1}}\overset{a}{F}|^{2}_{g^{a}(\cdot,t),h^{a}}+C|\overset{a}{\nabla^{k}}
      \overset{a}{F}|^{2}_{g^{a}(\cdot,t),h^{a}}  +C(n) \sum_{l=1}^{k-1}|\overset{a}{\nabla^{l}}[\overset{a}{R_M}+\overset{a}{R_{N}}\ast
      (\overset{a}{\nabla}\overset{a}{F})^{2}\\&\ \ +e^{\varphi^{a}}\overset{a}{R_M}
      +\overset{a}{\nabla^{2}}e^{\varphi^{a}}]|_{g^{a},h^{a}}|\overset{a}{\nabla^{k-l}}\overset{a}{F}|_{g^{a},h^{a}}
      |\overset{a}{\nabla^{k}}\overset{a}{F}|_{g^{a}(\cdot,t),h^{a}}.
        \end{split}
       \end{equation*}
     By the induction hypothesis, the local
derivative estimate of Shi, and (2.15), (2.19) and (2.20), it
follows
\begin{eqnarray*}
\sum_{l=0}^{k-1}|\overset{a}{\nabla^{l}}\overset{a}{R_{M}}|_{g
^{a}}|\overset{a}{\nabla^{k-l}}
\overset{a}{F}|_{g^{a}(\cdot,t),h^{a}}&\leqslant &\frac{C}{
t^{\frac{k-1}{2}}},\\
\sum_{l=0}^{k-1}|\overset{a}{\nabla^{l}}[\overset{a}{R_{N}}\ast(\overset{a}{\nabla}\overset{a}{F})^{2}]|_{g
^{a}}|\overset{a}{\nabla^{k-l}}
\overset{a}{F}|_{g^{a}(\cdot,t),h^{a}}&\leqslant &\frac{C}{
t^{\frac{k-1}{2}}}+C|\overset{a}{\nabla^{k}}
\overset{a}{F}|_{g^{a}(\cdot,t),h^{a}},\\
\sum_{l=0}^{k-1}|\overset{a}{\nabla^{l+2}}e^{\varphi^{a}}|_{g
^{a}}|\overset{a}{\nabla^{k-l}}
\overset{a}{F}|_{g^{a}(\cdot,t),h^{a}}&\leqslant &\frac{C}{
t^{\frac{k-2}{2}}},\\
\sum_{l=0}^{k-1}|\overset{a}{\nabla^{l}}\overset{a}{e^{\varphi^{a}}R_{M}}|_{g
^{a}}|\overset{a}{\nabla^{k-l}}
\overset{a}{F}|_{g^{a}(\cdot,t),h^{a}}&\leqslant &\frac{C}{
t^{\frac{k-1}{2}}}.\\
\end{eqnarray*}
This gives
\begin{equation*}      \begin{split}
     (\frac{\partial}{\partial
     t}-\overset{a}{\triangle}_{t})|\overset{a}{\nabla^{k}}
     \overset{a}{F}|^{2}_{g^{a}(\cdot,t),h^{a}}
     &\leqslant
      -2|\overset{a}{\nabla^{k+1}}\overset{a}{F}|^{2}_{g^{a}(\cdot,t),h^{a}}+C|\overset{a}{\nabla^{k}}
      \overset{a}{F}|^{2}_{g^{a}(\cdot,t),h^{a}}  +\frac{C}{t^{\frac{k-1}{2}}}
      |\overset{a}{\nabla^{k}}\overset{a}{F}|_{g^{a}(\cdot,t),h^{a}},\\
(\frac{\partial}{\partial
     t}-\overset{a}{\triangle}_{t})|\overset{a}{\nabla^{k}}
     \overset{a}{F}|_{g^{a}(\cdot,t),h^{a}}
     &\leqslant
      C|\overset{a}{\nabla^{k}}
      \overset{a}{F}|_{g^{a}(\cdot,t),h^{a}}
      +\frac{C}{t^{\frac{k-1}{2}}},
      \\
      (\frac{\partial}{\partial
     t}-\overset{a}{\triangle}_{t})|\overset{a}{\nabla^{k-1}}
     \overset{a}{F}|^{2}_{g^{a}(\cdot,t),h^{a}}
     &\leqslant
      -2|\overset{a}{\nabla^{k}}\overset{a}{F}|^{2}_{g^{a}(\cdot,t),h^{a}}+
      +\frac{C}{t^{k-\frac{5}{2}}}.
        \end{split}
       \end{equation*}
Let $\varepsilon=\frac{2(k-3)}{k-2}-1$, then $0\leq\varepsilon<1$
for $k\geq 4$. It is clear that
 $$(\frac{\partial}{\partial
     t}-\overset{a}{\triangle}_{t})|\overset{a}{\nabla^{k}}
     \overset{a}{F}|_{g^{a}(\cdot,t),h^{a}}^{1+\varepsilon}
     \leqslant
      C|\overset{a}{\nabla^{k}}
      \overset{a}{F}|_{g^{a}(\cdot,t),h^{a}}^{1+\varepsilon}  +\frac{C}{t^{\frac{k-1}{2}}}|\overset{a}{\nabla^{k}}
      \overset{a}{F}|_{g^{a}(\cdot,t),h^{a}}^{\varepsilon},$$
and
 \begin{eqnarray*}
 (\frac{\partial}{\partial
     t}-\overset{a}{\triangle}_{t})(|\overset{a}{\nabla^{k}}
     \overset{a}{F}|_{g^{a}(\cdot,t),h^{a}}^{1+\varepsilon}+|\overset{a}{\nabla^{k-1}}
     \overset{a}{F}|^{2}_{g^{a}(\cdot,t),h^{a}})
     &\leqslant&
     -2|\overset{a}{\nabla^{k}}
     \overset{a}{F}|^{2}_{g^{a}(\cdot,t),h^{a}}
      C|\overset{a}{\nabla^{k}}
      \overset{a}{F}|_{g^{a}(\cdot,t),h^{a}}^{1+\varepsilon}  \\& &+\frac{C}{t^{\frac{k-1}{2}}}|\overset{a}{\nabla^{k}}
      \overset{a}{F}|_{g^{a}(\cdot,t),h^{a}}^{\varepsilon}+\frac{C}{t^{k-\frac{5}{2}}},
\end{eqnarray*}
on
$B_{g^{a}(\cdot,0)}(x_{0},(\frac{1}{2}+\frac{1}{2^{k}})\delta)\times[0,\bar{T}]$.

Let $$v=t^{k-3}(|\overset{a}{\nabla^{k}}
     \overset{a}{F}|_{g^{a}(\cdot,t),h^{a}}^{1+\varepsilon}+|\overset{a}{\nabla^{k-1}}
     \overset{a}{F}|^{2}_{g^{a}(\cdot,t),h^{a}}).$$
 Then we have
\begin{eqnarray*}
(\frac{\partial}{\partial t}-\overset{a}{\triangle})v&\leqslant&
(k-3)\frac{v}{t}+ t^{k-3}(-|\overset{a}{\nabla^{k}}
     \overset{a}{F}|^{2}_{g^{a}(\cdot,t),h^{a}}
      +\frac{C}{t^{\frac{k-1}{2}}}|\overset{a}{\nabla^{k}}
      \overset{a}{F}|_{g^{a}(\cdot,t),h^{a}}^{\varepsilon}+\frac{C}{t^{k-\frac{5}{2}}}\\ &\leqslant&
      -\frac{1}{t}\{v^{\frac{2}{1+\varepsilon}}
      -C\sqrt{t}v^{\frac{\varepsilon}{1+\varepsilon}}-C\sqrt{t}\}\\
      &\leqslant&-\frac{1}{t}\{v^{\frac{2}{1+\varepsilon}}
      -C\}
\end{eqnarray*}
on
$B_{g^{a}(\cdot,0)}(x_{0},(\frac{1}{2}+\frac{1}{2^{k}})\delta)\times[0,\bar{T}]$.
Similarly, at the maximum point of
$\xi_{k}(\frac{d_{g^{a}(\cdot,0)}(x_{0},\cdot)}{\delta})v$, we
have
\begin{equation*}
\begin{split}0&\leqslant-\frac{1}{t}\{\xi_{k}v^{\frac{2}{1+\varepsilon}}-C\xi_{k}\}
+C(\frac{|\xi_{k}^{\prime}|^{2}}{\xi_{k}}+|\xi_{k}^{\prime\prime}|)v\\
&\leqslant-\frac{1}{t}\{\xi_{k}v^{\frac{2}{1+\varepsilon}}-C\xi_{k}^{\frac{1+\varepsilon}{2}}v-C\}\\
&\leqslant-\frac{1}{t}\{\frac{1}{2}\xi_{k}v^{\frac{2}{1+\varepsilon}}-C\}\\
&\leqslant-\frac{1}{t}\{\frac{1}{2}(\xi_{k}v)^{\frac{2}{1+\varepsilon}}-C\}.
\end{split}
\end{equation*}
since $\frac{2}{1+\varepsilon}>1$. So
 we proved the $k$-th order estimate
 $$ |\overset{a}{\nabla^{k}} \overset{a}{F}|_{g^{a}(\cdot,t),h^{a}}\leq C t^{-\frac{k-2}{2}}
           $$
    on $B_{g^{a}(\cdot,0)}(x_{0},(\frac{1}{2}+\frac{1}{2^{k+1}})\delta)\times [0,\bar{T}]$.
This completes the proof of the lemma.

$\hfill\#$

Now we are ready to prove Theorem 2.7.

\vskip 0.5cm $\underline{\mbox{\textbf{Proof of Theorem 2.7.}}}$
\vskip 0.3cm

Since $D_j\supseteq B_{g^{a}(\cdot,0)}(P,j+1),$ by choosing
$\delta=1$ and $\bar{T}=T_4$ in Lemma 2.11 and Lemma 2.12, we get
a convergent subsequence of $\overset{a}{F^{j}}$ (as
$j\rightarrow\infty$) on $B_{g^{a}(\cdot,0)}(P,j)\times [0,T_4]$.
Denote the limit by $\overset{a}{F}$ (as $j\rightarrow\infty$).
Then $\overset{a}{F}$ is the desired solution of
$(2.3)_{a}^{\prime}$ with estimates (2.22).

Finally we prove a uniqueness theorem for the solutions of
$(2.3)_{a}^{\prime}$
 with estimates (2.22).
\begin{lem}Let $\overset{a}{F}$ and $\overset{a}{\bar{F}}$ be two solutions
of the intial problem $(2.3)_{a}^{\prime}$ on $[0,\bar{T}]$,
$\bar{T}\leqslant T$, with estimates (2.22). Then
$\overset{a}{F}=\overset{a}{\bar{F}}$ on $[0,\bar{T}]$.
\end{lem}
$\underline{\mbox{\textbf{Proof }}}$  Set
$\psi(y_{1},y_{2})=\frac{1}{2}d^2_{(N^{n},h^{a})}(y_{1},y_{2})$
  and $\rho(x,t)=\psi(\overset{a}{F}(x,t),\overset{a}{\bar{F}}(x,t))$. Then $\psi(x,t)$ is smooth
  when $\psi<\frac{1}{2}i_{a}^{2}$.
   Now by the same calculation as in
Lemma 2.8, we have:
$$ (\frac{\partial}{\partial
     t}-\overset{a}{\triangle_{t}})\rho=-Hess(\psi)(X_{i},X_{j})(g^{a})^{ij}
$$ where the vector fields $X_{i}$, $i=1, 2, \cdots, n$,
in local coordinates $(y_{1}^{\alpha},y_{2}^{\beta})$ on
$N^{n}\times N^{n}$ are defined as follows
$$
X_{i}=\frac{\partial \overset{a}{F^{\alpha}}}{\partial
x^{i}}\frac{\partial}{\partial y_{1}^{\alpha}}+\frac{\partial
\overset{a}{\bar{F}^{\beta}}}{\partial
x^{i}}\frac{\partial}{\partial y_{2}^{\beta}}.
$$
By the estimates (2.22), we know that there is a constant
$0<\bar{T^{\prime}}\leqslant\bar{T}$ such that there holds
$$\rho<\min\{\frac{i_{a}^{2}}{8},\frac{\pi^{2}}{8\bar{C_0}}\}.$$
on $M^n\times [0,\bar{T^{\prime}}]$.

 Similarly as in the proof of Lemma 2.8. By using the computation of $Hess(\psi)$ in \cite{ScY1} (the formula (16) in
  \cite{ScY1}), for any $(u,v)\in
D=\{(u,v): (u,v)\in N^{n}\times N^{n} \mbox{ with }
d_{N^{n}}(u,v)<\min\{\frac{i_{a}}{2},\frac{\pi}{2\sqrt{\bar{C_0}}}\}\}$,
and any $X \in T_{(u,v)}D$,
\begin{eqnarray*}
Hess(\psi)(X,X)&\geqslant &-\int_{0}^{r}t\langle
R(e_{1},V)V,e_{1}\rangle-\int_{0}^{r}t\langle
R(\bar{e}_{1},V)V,\bar{e}_{1}\rangle
\end{eqnarray*}
where $V$ is a Jacobi field on geodesic $\sigma$ (connecting
$(v,v)$ to $(u,v)$) and $\bar{\sigma}$ (connecting $(u,u)$ to
$(u,v)$) with $X$ as the boundary values. Since
$|\overset{a}{\nabla} F|_{g^{a},h^{a}}$ and
$|\overset{a}{\nabla}\bar{F}|_{g^{a},h^{a}}$ are bounded,  we know
from above formula that
$$Hess(\psi)(X_{i},X_{j})(g^{a})^{ij}\geqslant -C\rho$$
on $M^{n} \times[0,\bar{T^{\prime}}]$.  Thus we have
$$ (\frac{\partial}{\partial
     t}-\overset{a}{\triangle_{t}})\rho \leqslant C\rho
$$
on $M^{n} \times[0,\bar{T^{\prime}}]$. By the maximum principle,
it follows that $\rho=0$ on $M^{n} \times[0,\bar{T^{\prime}}]$.
Then the lemma follows by continuity method.

 $\hfill\#$

\subsubsection{Proof of theorem 2.6 and Theorem 2.1}
\vskip 0.5cm $\underline{\mbox{\textbf{Proof of Theorem
2.6.}}}$\vskip 0.3cm

Let us check the initial data. Now $f=identity$, so
\begin{equation*} \tag{2.36}
\begin{split}|\overset{a}{\nabla}f|^{2}_{g^{a}(\cdot,0),h^{a}}&=
g^{ij}(\cdot,0)g_{ij}(\cdot,T)\\
&\leqslant ne^{2nk_0T}
\end{split}
\end{equation*}
\begin{equation*} \tag{2.37}
\begin{split}|\overset{a}{\nabla^{2}}f|^{2}_{g^{a}(\cdot,0),h^{a}}
&=|\overset{a}{\Gamma^{k}_{ij}}(\cdot,0)-\overset{a}{\Gamma^{k}_{ij}}(\cdot,T)|_{g^{a}(\cdot,0),h^{a}}
\\
&\leqslant
C(n,k_0,T)\int_{0}^{T}e^{\varphi^{a}}(|\overset{a}{\nabla}\overset{a}{R}_{M}|_{g^{a}(\cdot,t)}+
|\overset{a}{R}_{M}\ast\overset{a}{\nabla}\varphi^{a}|_{g^{a}(\cdot,t)})\\&
\ \ \
+|\overset{a}{\nabla}\varphi^{a}|_{g^{a}(\cdot,t)}|\overset{a}{\nabla^{2}}\varphi^{a}|_{g^{a}(\cdot,t)}+
|\overset{a}{\nabla^{3}}\varphi^{a}|_{g^{a}(\cdot,t)})dt\\
&\leqslant C(n,k_0,T)\int_{0}^{T}\frac{1}{\sqrt{t}}+|\log t|dt\\
& \leqslant C(n,k_0,T).
\end{split}
\end{equation*}
By applying Theorem 2.7, we know that there is $\delta_0>0$ such
that $(2.3)_a$ has a smooth solution $\overset{a}{F}$ on
$M^{n}\times[0,\delta_0]$ with estimates (2.22). In views of Lemma
2.12 and Lemma 2.13, in order to prove Theorem 2.6, we only need
to bound
$|\overset{a}{\nabla}\overset{a}{F}|^{2}_{g^{a}(\cdot,t),h^{a}}$
 uniformly on a uniformly interval $[0,T_1]$ with $T_1$ independent
 of $a$. To this end, let
\begin{equation*}
\begin{split}
 \tilde{T}=\sup\{\tilde{T_0}\mid &\tilde{T_0}\leqslant T, (2.3)_a
 \text{ has a smooth solution on } M^{n}\times [0,\tilde{T_0}]\\
 &\text{with}
 \sup_{M^{n}\times [0,\tilde{T_0}]}|\overset{a}{\nabla}\overset{a}{F}|^{2}_{g^{a}(\cdot,t),h^{a}}
 <\infty\},
\end{split}
\end{equation*}
We will estimate $\tilde{T}$ from below.

We come back to the equation (2.30) of
$|\overset{a}{\nabla}\overset{a}{F}|^{2}_{g^{a}(\cdot,t),h^{a}}$,
where there holds
\begin{equation*}     \begin{split}
     (\frac{\partial}{\partial
     t}-\overset{a}{\triangle}_{t})|\overset{a}{\nabla}
     \overset{a}{F}|^{2}_{g^{a}(\cdot,t),h^{a}}
     \leqslant
      -2|\overset{a}{\nabla^{2}}\overset{a}{F}|^{2}_{g^{a}(\cdot,t),h^{a}}+C_1(n,k_0,T)|\overset{a}{\nabla}
      \overset{a}{F}|^{2}_{g^{a}(\cdot,t),h^{a}}+C_2(n,k_0,T)
      |\overset{a}{\nabla}\overset{a}{F}|^{4}_{g^{a}(\cdot,t),h^{a}}
        \end{split}
       \end{equation*}
on $ M^{n}\times [0,\tilde{T}]$.  We remark that $\overset{a}{F}$
is defined on a complete manifold with bounded curvature and
$\sup_{M^{n}\times
[0,\tilde{T_0}]}|\overset{a}{\nabla}\overset{a}{F}|^{2}_{g^{a}(\cdot,t),h^{a}}
 <\infty$, for each $\tilde{T_0}<\tilde{T}$. So by applying the
 maximum principle on complete manifolds, we have
 $$\frac{d^{+}}{d
 t}(\sup_{M^{n}}|\overset{a}{\nabla}\overset{a}{F}|^{2}_{g^{a}(\cdot,t),h^{a}})\leqslant C_1(n,k_0,T)
 \sup_{M^{n}}|\overset{a}{\nabla}\overset{a}{F}|^{2}_{g^{a}(\cdot,t),h^{a}}+
 C_2(n,k_0,T)\sup_{M^{n}}|\overset{a}{\nabla}\overset{a}{F}|^{4}_{g^{a}(\cdot,t),h^{a}}$$
 where $\frac{d^{+}}{d
 t}$ is the upper right derivative defined by
 $$
 \frac{d^{+}}{d
 t} u=\limsup_{\triangle t\searrow 0}\frac{u(t+\triangle t)-u(t)}{\triangle
 t}.
 $$
 By combining with (2.36), we have
$$\sup_{M^{n}\times
[0,\tilde{T_0}]}|\overset{a}{\nabla}\overset{a}{F}|^{2}_{g^{a}(\cdot,t),h^{a}}\leqslant
2n e^{2nk_0T},$$
provided $\tilde{T_0}\leqslant \min\{T,
\frac{\log 2}{C_1(n,k_0,T)+2ne^{2nk_0T}C_2(n,k_0,T)}\}$.

By Lemma 2.12 and Lemma 2.13 and Theorem 2.7,  the solution
$\overset{a}{F}$ exists smoothly until
$|\overset{a}{\nabla}\overset{a}{F}|^{2}_{g^{a}(\cdot,t),h^{a}}$
blows up, so we know $\tilde{T}\geqslant \min\{T, \frac{\log
2}{C_1(n,k_0,T)+2ne^{2nk_0T}C_2(n,k_0,T)}\}$. By choosing
$T_1=\min\{T, \frac{\log
2}{C_1(n,k_0,T)+2ne^{2nk_0T}C_2(n,k_0,T)}\}$, Theorem 2.6 follows.

$\hfill\#$

 \vskip 0.5cm $\underline{\mbox{\textbf{Proof of Theorem 2.1.}}}$
 \vskip 0.3cm

Note that $\varphi^{a}=0$ on $B_{g(\cdot,T)}(P,a)$, and
$g^{a}_{ij}(x,t)=e^{\varphi^{a}}g_{ij}(x,t)$,
$h^{a}_{\alpha\beta}(y)=e^{\varphi^{a}}h_{\alpha\beta}$. It
follows that
$$
g^{a}_{ij}(x,t)=g_{ij}(x,t)   \ \ \  \text{on} \
B_{g(\cdot,T)}(P,a),
$$
$$
h^{a}_{\alpha\beta}(y)=h_{\alpha\beta}(y)   \ \ \  \text{on} \
B_{g(\cdot,T)}(P,a).
$$
By Theorem 2.6 and estimates (2.21) and letting $a\rightarrow
\infty$, the solutions $\overset{a}{F}$ of $(2.3)_{a}$ on
$M^{n}\times[0,T_1]$ have a convergent subsequence so that the
limit is a
 solution of (2.3) with the estimates (2.4).

 $\hfill\#$
\vskip 0.8cm
\section{The uniqueness of the Ricci flow }
\vskip 0.5cm
\subsection{Preliminary estimates for the Ricci-De Turck flow}
  Let $F(x,t)$ be a solution to (2.3) in Theorem 2.1 on
  $M^{n}\times[0,T_0]$. Let
  $\tilde{g}_{ij}(x,t)=h_{\alpha\beta}(F(x,t))\frac{\partial F^{\alpha}}{\partial x^{i}}
  \frac{\partial F^{\beta}}{\partial x^{j}}$ be the one-parameter
  family of pulled back metrics $F^{*}{h}$. We will estimate $g_{ij}(x,t)$ in terms of
  $\tilde{g}_{ij}(x,t)$.
 \begin{prop}
 There exists a constant
 $0<T_5\leq T_0$ depending only on $k_0$ and $T$ such that for all $(x,t)\in M^{n}\times
 [0,T_5]$, we have
 \begin{equation*}\tag{3.1}
 \begin{split}
 \frac{1}{C(n,k_0,T)}\tilde{g}_{ij}(x,t)\leq g_{ij}(x,t)&\leq
 C(n,k_0,T)\tilde{g}_{ij}(x,t)\\
 |\tilde{\nabla}^{k}g|_{\tilde{g}}\leq \frac{C(n,k_0,T,k)}{t^{\frac{k-1}{2}}}
 \end{split}
 \end{equation*}
 for $k=1,2,\cdots$
 \end{prop}
 Proof.  We first consider the zero-order estimate of
 $g_{ij}(x,t)$. The estimate $|\nabla F|^{2}=\tilde{g}_{ij}g^{ij}\leq C$
in (2.4) implies $\tilde{g}_{ij}(x,t)\leq Cg_{ij}(x,t)$. For the
reverse inequality, we compute the equation of
$\tilde{g}_{ij}(x,t)$:
\begin{eqnarray*}
\frac{\partial}{\partial t}\tilde{g}_{ij}&=&\triangle
\tilde{g}_{ij}-2R_{ik}F^{\alpha}_{l}F^{\beta}_{j}h_{\alpha\beta}g^{kl}
+2R_{\alpha\beta\gamma\delta}F^{\alpha}_{i}F^{\beta}_{k}
F^{\gamma}_{j}F^{\delta}_{l}g^{kl}-2h_{\alpha\beta}F^{\alpha}_{k,i}F^{\beta}_{l,j}g^{kl}\\
&\geq&\triangle \tilde{g}_{ij}-
2R_{ik}\tilde{g}_{jl}g^{kl}-2k_0|\nabla
F|^{2}g_{ij}-2|\nabla^{2}F|^{2}g_{ij}\\
&\geq&\triangle
\tilde{g}_{ij}-2R_{ik}\tilde{g}_{jl}g^{kl}-C(n,k_0,T) g_{ij},
\end{eqnarray*}
by (2.4). Combining this with the Ricci flow equation gives
\begin{eqnarray*}
(\frac{\partial}{\partial t}-\triangle)(\tilde{g}_{ij}+
C(n,k_0,T)t
 g_{ij}-\frac{1}{2ne^{2nk_0T}}g_{ij})\geqslant-2R_{ik}(\tilde{g}_{lj}+
C(n,k_0,T)t
 g_{lj}-\frac{1}{2ne^{2nk_0T}}g_{lj})g^{kl}.
\end{eqnarray*}
Note that at $t=0$,
$$(\tilde{g}_{ij}+
C(n,k_0,T)t
 g_{ij}-\frac{1}{2ne^{2nk_0T}}g_{ij})\mid_{t=0}=g_{ij}(\cdot,T)-\frac{1}{2ne^{2nk_0T}}g_{ij}(\cdot,0)>0.
$$
By applying  the maximum principle to  above equation, we obtain
$$
\tilde{g}_{ij}+ C(n,k_0,T)t
 g_{ij}-\frac{1}{2ne^{2nk_0T}}g_{ij}>0
$$
on $M^{n}\times[0,T_0]$. Let
$T_5=\min\{T_0,\frac{1}{4ne^{2nk_0T}C(n,k_0,T)}\}$. Then we have
$$\tilde{g}_{ij}\geq \frac{1}{4ne^{2nk_0T}} g_{ij}, \ \ \ \text{on}\ \ M^{n}\times[0,T_5].$$ This gives the zero-order estimate
of $g_{ij}(x,t)$.

 For the first order derivative of
$g_{ij}$, we compute
\begin{eqnarray*}
\tilde{\nabla}_{k}g_{ij}=(\tilde{\nabla}_{k}-\nabla_{k})g_{ij}
=(\Gamma^{l}_{ki}-\tilde{\Gamma}^{l}_{ki})g_{lj}+(\Gamma^{l}_{kj}-\tilde{\Gamma}^{l}_{kj})g_{li}
\end{eqnarray*}
and
\begin{eqnarray*}
|(\Gamma^{l}_{ki}-\tilde{\Gamma}^{l}_{ki})|^{2}_{\tilde{g}}&=
&|(\Gamma^{p}_{ki}-\tilde{\Gamma}^{p}_{ki})\tilde{g}_{lp}|^{2}_{\tilde{g}}\\
&=&|\nabla_k\nabla_iF^{\alpha}\frac{\partial F^{\beta}}{\partial
x^{l}}h_{\alpha\beta}|_{\tilde{g}}\\
&\leqslant& C(n,k_0,T)|\nabla_k\nabla_iF^{\alpha}\frac{\partial
F^{\beta}}{\partial
x^{l}}h_{\alpha\beta}|_{g}\\
&\leqslant& C(n,k_0,T)|\nabla^{2}F|_{g,h}|\nabla F|_{g,h}\\
&\leqslant& C(n,k_0,T).
\end{eqnarray*}
This gives the first order estimate.

For higher order estimates, we prove it by induction. Suppose we
have showed
$$
|\tilde{\nabla^{l}}g|_{\tilde{g}}\leqslant\frac{C}{t^{\frac{l-1}{2}}}
\ \ \ \ \text{for}\ \ l=1,2,\cdots,k-1,
$$
$$
|\tilde{\nabla^{l}}(\Gamma-\tilde{\Gamma})|_{\tilde{g}}\leqslant\frac{C}{t^{\frac{l}{2}}}
\ \ \ \ \text{for}\ \ l=0,1,\cdots,k-2.
$$
Since by induction
 \begin{eqnarray*}
|\tilde{\nabla}^{k-1}(\Gamma-\tilde{\Gamma})|_{\tilde{g}}&=&
|\tilde{\nabla}^{k-1}[(\Gamma-\tilde{\Gamma})\ast\tilde{g}]|_{\tilde{g}}\\
 &=&|\sum_{j=0}^{k-1}\nabla^{k-1-j}[(\Gamma-\tilde{\Gamma})\ast\tilde{g}]\ast\sum_{i_1+1+\cdots+i_q+1=j
}\tilde{\nabla}^{i_1}(\Gamma-\tilde{\Gamma})\ast\cdots\ast\tilde{\nabla}^{i_q}(\Gamma-\tilde{\Gamma})|_{\tilde{g}}\\
&\leqslant& C(n,k_0,T)
\sum_{j=0}^{k-1}|\nabla^{k-1-j}(\nabla^{2}F\ast\nabla F)|_{g,h}\\
& & \times\sum_{i_1+1+\cdots+i_q+1=j}|\tilde
{\nabla}^{i_1}(\Gamma-\tilde{\Gamma})|_{\tilde{g}}\cdots |\tilde
{\nabla}^{i_q}(\Gamma-\tilde{\Gamma})|_{\tilde{g}}\\
&\leqslant&
C(n,k_0,T,k)(\frac{1}{t^{\frac{k-1-j}{2}}}\frac{1}{t^{\frac{j-2}{2}}}+\frac{1}{t^{\frac{k-1}{2}}})\\
&\leqslant& \frac{C(n,k_0,T,k)}{t^{\frac{k-1}{2}}}
\end{eqnarray*}
 and
\begin{eqnarray*}
\tilde{\nabla}^{k}g&=&\tilde{\nabla}^{k-1}((\Gamma-\tilde{\Gamma})\ast
g
)\\&=&\sum_{i=0}^{k-1}\tilde{\nabla^{i}}(\Gamma-\tilde{\Gamma})\ast\tilde{\nabla}^{k-1-i}g,
\end{eqnarray*} then  we have
$$
|\tilde{\nabla}^{k}g|_{\tilde{g}}\leqslant
\frac{C}{t^{\frac{k-1}{2}}}.
$$
This completes the induction argument and the proposition is
proved.

$\hfill\#$

\begin{prop}
Let $F(x,t)$ be the solution of (2.3) in Theorem 2.1. Then
$F(\cdot,t)$ are diffeomorphisms for all $t\in [0,T_5]$; moreover,
there exists a constant $C(n,k_0,T)>0$ depending only on $n$,
$k_0$ and $T$ such that
$$
d_{h}(F(x_1,t),F(x_2,t))\geqslant e^{-C(n,k_0,T)}d_{h}(x_1,x_2)
$$
for all $x_1,x_2\in M^{n}$, $t\in[0,T_{5}]$.
\end{prop}
Proof. Note that
$$
\frac{1}{C}\tilde{g}_{ij}(x,t)\leq g_{ij}(x,t)\leq
C\tilde{g}_{ij}(x,t)
$$ implies that $F$ are local diffeomorphisms. So we only need to prove that $F(\cdot,t)$ is injective.
Suppose not. Then there exist two points  $x_1\neq x_2$, such that
$ F(x_1,t)= F(x_2,t)$, for some $t_0\in(0,T_5]$. Assume $t_0>0$ be
the first time so that $ F(x_1,t)= F(x_2,t)$. Choose small
$\delta>0$, such that there exist a neighborhood $\tilde{O}$ of
$F(x_1,t_0)$ and a neighborhood $O$ of $x_1$ such that
$F^{-1}(\cdot,t)$ is a diffeomorphism from $\tilde{O}$ to $O$ for
all $t\in[t_0-\delta,t_0]$, moreover, letting $\tilde{\gamma_t}$
be a shortest geodesic( parametrized by arc length)
 on the target $(N^{n},h_{\alpha\beta})$ connecting  $F(x_1,t)$ and $F(x_2,t)$,
 we require $\tilde{\gamma}\in \tilde{O}$ for $t\in[t_0-\delta,t_0]$. We
 compute
 \begin{eqnarray*}\frac{\partial}{\partial t}d_{h}(F(x_1,t),F(x_2,t))&=&\langle
 V,{\tilde{\gamma}}^{\prime}(l)\rangle_{h}-\langle
 V,{\tilde{\gamma}}^{\prime}(0)\rangle_{h}
  \end{eqnarray*}
  where $\tilde{\gamma}(0)=F(x_1,t)$ ,
  $\tilde{\gamma}(l)=F(x_2,t)$, and $V^{\alpha}=\triangle
  F^{\alpha}$.
 Now we pull back everything by $F^{-1}$ to $O$,
 \begin{eqnarray*}\frac{\partial}{\partial t}d_{h}(F(x_1,t),F(x_2,t))&=&\langle
 P_{-\tilde{\gamma}}V-V,{\gamma}^{\prime}(0)\rangle_{F^{*}h}\\
 &\geq&- \sup_{x\in F^{-1}\tilde{\gamma}}|\tilde{\nabla}V|(x,t) d_{h}(F(x_1,t),F(x_2,t))
 \end{eqnarray*}
 where $P_{\tilde{\gamma}}$ is the parallel translation along $F^{-1}\tilde{\gamma}$
 using the metric $F^{*}h$. By (2.4),
 \begin{eqnarray*}
 |\tilde{\nabla}_{k} V^{l}|_{\tilde{g}}&=&|\tilde{\nabla}_{k}(\triangle
 F^{\alpha}\frac{\partial F^{\alpha}}{\partial x^{l} }h_{\alpha\beta})|_{\tilde{g}}\\&\leqslant&
 |\nabla_{k}(\triangle
 F^{\alpha}\frac{\partial F^{\alpha}}{\partial x^{l} }h_{\alpha\beta})|_{\tilde{g}}+
 C|\Gamma-\tilde{\Gamma}|_{\tilde{g}}|\nabla^{2}F|_{g,h}|\nabla F|_{g,h}\\
 &\leqslant&
 C(n,k_0,T)(|\nabla^{3}F|_{g,h}|\nabla F|_{g,h}+|\nabla^{2}F|^{2}_{g,h}+|\nabla^{2}F|^{2}_{g,h}|\nabla F|_{g,h})\\
 &\leqslant& \frac{C(n,k_0,T)}{\sqrt{t}}.
 \end{eqnarray*}
 It follows that  we have
 $$d_{h}(F(x_1,t),F(x_2,t))\leqslant e^{C(\sqrt{t_0}-\sqrt{t_0-\delta})
 }d_{h}(F(x_1,t_0),F(x_2,t_0))=0,$$
 for $t\in [t_0-\delta,t_0]$, which
 contradicts with the choice of $t_0$. So $F(\cdot,t)$ are diffeomorphisms.

 By choosing $\tilde{O}=N^{n}$, $O=M^{n}$, the above computation also gives
$$
d_{h}(F(x_1,t),F(x_2,t))\geqslant e^{-C(n,k_0,T)}d_{h}(x_1,x_2).
$$ The proof of
 the proposition is completed.

 $\hfill\#$

\subsection{Ricci-De Turck flow}

 \ \ From the previous section, we know that the harmonic map flow coupled
with Ricci flow (2.3)with identity as initial data has a short
time solution $F(x,t)$ on $M^{n}\times[0,T_5]$, which remains
being a diffeomorphism with good estimates (2.4). Let
${(F^{-1})}^{*}g$ be one-parameter  family of pulled back metrics
on the target $(N^{n},h_{\alpha\beta})$. Denote
$g_{\alpha\beta}(y,t)$. Then $g_{\alpha\beta}(y,t)$ satisfies the
so called Ricci-De Turck flow:
$$
\frac{\partial}{\partial
t}g_{\alpha\beta}(y,t)=-2R_{\alpha\beta}(y,t)+\nabla_{\alpha}
V_{\beta}+\nabla_{\beta} V_{\alpha} \eqno(3.3)
$$
where
$V^{\alpha}=g^{\beta\gamma}(\Gamma^{\alpha}_{\beta\gamma}(g)-{\Gamma}^{\alpha}_{\beta\gamma}(h))$,
$\Gamma^{\alpha}_{\beta\gamma}(g)$ and
${\Gamma}^{\alpha}_{\beta\gamma}(h)$ are the Christoffel symbols
of the metrics $g_{\alpha\beta}(y,t)$ and $h_{\alpha\beta}(y)$
respectively.

 By (3.1) of Proposition 3.1, we already have the following estimates for $g_{\alpha\beta}(y,t)$

  $$\frac{1}{C(n,k_{0},T)}h_{\alpha\beta}(y)\leqq g_{\alpha\beta}(y,t) \leqq C(n,k_{0},T) h_{\alpha\beta}(y)$$
 $$|{\nabla^{k}_{h}}g|_{h}\leqq
 \frac{C(n,k_{0},T,k)}{t^{\frac{k-1}{2}}}.\eqno(3.4)$$
on $N^{n}\times[0,T_5]$.

  Let $g_{ij}(x,t)$ and
$\bar{g}_{ij}(x,t)$ be two solutions to the Ricci flow with
bounded curvature and with the same initial value as assumed in
Theorem 1.1. We solve the corresponding harmonic map flow with
same target $(N^{n},h_{\alpha\beta})=(M^{n},g_{ij}(\cdot,T))$ by
\begin{equation*} \tag{3.5}
\left\{
\begin{split}
 \quad \frac{\partial}{\partial t}F(x,t)&=\triangle F(x,t),  \\
  F(\cdot,0)&= identity,
  \end{split}
  \right.
\end{equation*} and
\begin{equation*} \tag{3.6}
\left\{
\begin{split}
 \quad \frac{\partial}{\partial t}\bar{F}(x,t)&=\bar{\triangle}\bar{F}(x,t),  \\
  \bar{F}(\cdot,0)&=identity,
\end{split}
  \right.
\end{equation*}
respectively. Then we obtained two solutions $F(x,t)$ and
$\bar{F}(x,t)$ on $M^{n}\times[0,T_5]$.  It is clear that
$\bar{F}(x,t)$ still satisfies (2.4), Proposition 3.1 and
Proposition 3.2.  Let
$\bar{g}_{\alpha\beta}(y,t)={(\bar{F}^{-1})}^{*}\bar{g}(y,t)$,
then $\bar{g}_{\alpha\beta}(y,t)$ still satisfies (3.4). Now we
have two solutions $g_{\alpha\beta}(y,t)$ and
$\bar{g}_{\alpha\beta}(y,t)$ to the Ricci De-Turck flow with same
initial data and with good estimates (3.4).
\begin{prop} There holds
   $$g_{\alpha\beta}(y,t)=\bar{g}_{\alpha\beta}(y,t)$$  on
$N^{n}\times [0,T_5]$.
\end{prop}
 Proof. We can write the Ricci-De Turck flow (3.3) by
 using the fixed metric $h_{\alpha\beta}(y)$ in the following form (see \cite{Sh1}):
 \begin{equation*} \tag{3.7}
\begin{split}
 \quad \frac{\partial}{\partial t}g_{\alpha\beta}=&g^{\gamma\delta}\tilde{\nabla}_{\gamma}\tilde{\nabla}_{\delta}g_{\alpha\beta}-
 g^{\gamma\delta}g_{\alpha\xi}\tilde{g}^{\xi\eta}\tilde{R}_{\beta\gamma\eta\delta}
  -g^{\gamma\delta}g_{\beta\xi}\tilde{g}^{\xi\eta}\tilde{R}_{\alpha\gamma\eta\delta}+\frac{1}{2}g^{\gamma\delta}g^{\xi\eta}
  (\tilde{\nabla}_{\alpha}g_{\xi\gamma}\tilde{\nabla}_{\beta}g_{\eta\delta}\\&+
  2 \tilde{\nabla}_{\gamma}g_{\beta\xi}\tilde{\nabla}_{\eta}g_{\alpha\delta}-2\tilde{\nabla}_
  {\gamma}g_{\beta\xi}\tilde{\nabla}_{\delta}g_{\alpha\eta}-2\tilde{\nabla}_{\beta}
  g_{\xi\gamma}\tilde{\nabla}_{\delta}g_{\alpha\eta}
   -2 \tilde{\nabla}_{\alpha}g_{\xi\gamma}\tilde{\nabla}_{\delta}g_{\beta\eta})
  \end{split}
 \end{equation*}
  where $\tilde{g}_{\alpha\beta}=h_{\alpha\beta},$
  $\tilde{\nabla}$ and $\tilde{R}$ are the covariant derivative
  and the curvature of $\tilde{g}_{\alpha\beta}$.
  Note that $\bar{g}_{\alpha\beta}$  also satisfies (3.7), then the difference
  $g_{\alpha\beta}-\bar{g}_{\alpha\beta}$ satisfies the following
  equation:
\begin{equation*} \tag{3.8}
\begin{split}
  \frac{\partial}{\partial
  t}(g-\bar{g})=&g^{\gamma\delta}\tilde{\nabla}_{\gamma}\tilde{\nabla}_{\delta}(g-\bar{g})
  +g^{-1}\ast\bar{g}^{-1}\ast\tilde{\nabla}^{2}\bar{g}\ast(\bar{g}-g)\\
  &+ \bar{g}^{-1}\ast\tilde{g}^{-1}\ast\tilde{Rm}\ast(g-\bar{g})+g^{-1}\ast\bar{g}^{-1}\ast{g}\ast\tilde{g}^{-1}\ast\tilde{Rm}\ast(g-\bar{g})\\
  & +g^{-1}\ast
  g^{-1}\ast\bar{g}^{-1}\ast\tilde{\nabla}g\ast\tilde{\nabla}g\ast(g-\bar{g})+g^{-1}\ast
  \bar{g}^{-1}\ast\bar{g}^{-1}\ast\tilde{\nabla}g\ast\tilde{\nabla}g\ast(g-\bar{g})\\&
   +
  \bar{g}^{-1}\ast\bar{g}^{-1}\ast\tilde{\nabla}g\ast\tilde{\nabla}(g-\bar{g})+\bar{g}^{-1}\ast\bar{g}^{-1}\ast\tilde{\nabla}\bar{g}\ast\tilde{\nabla}(g-\bar{g})
  \end{split}
  \end{equation*}
   since $g^{\alpha\beta}-\bar{g}^{\alpha\beta}=g^{\alpha\xi}\bar{g}^{\eta\beta}(\bar{g}_{\eta\xi}-g_{\eta\xi})$.
   Let $$|g-\bar{g}|^{2}=\tilde{g}^{\alpha\gamma}\tilde{g}^{\beta\delta}
 (g_{\alpha\beta}-\bar{g}_{\alpha\beta})(g_{\gamma\delta}-\bar{g}_{\gamma\delta}).$$  It follows from (3.8) that:
 \begin{eqnarray*}
 (\frac{\partial}{\partial
 t}-g^{\gamma\delta}\tilde{\nabla}_{\gamma}\tilde{\nabla}_{\delta})|g-\bar{g}|^{2}
 &\leqslant&-2g^{\xi\eta}\tilde{g}^{\alpha\gamma}\tilde{g}^{\beta\delta}
 (\tilde{\nabla}_{\xi}g_{\alpha\beta}-\tilde{\nabla}_{\xi}\bar{g}_{\alpha\beta})
 (\tilde{\nabla}_{\eta}g_{\gamma\delta}-\tilde{\nabla}_{\eta}\bar{g}_{\gamma\delta})\\&
 & +100[|\tilde{Rm}|(1+|g||g^{-1}|)|\bar{g}^{-1}|
 +|\tilde{\nabla}^{2}\bar{g}||\bar{g}^{-1}||g^{-1}|\\& & \ \ \ \ \ +|\tilde{\nabla}g|^{2}
 (|\bar{g}^{-1}|^{2}|g^{-1}|+|\bar{g}^{-1}||g^{-1}\|^{2})]|g-\bar{g}|^{2}\\&
 & +100|\bar{g}^{-1}|^{2}(|\tilde{\nabla}g|+|\tilde{\nabla}\bar{g}|)|\tilde{\nabla}(g-\bar{g})||g-\bar{g}|
 \end{eqnarray*}
 where all the norms are computed with the metric $\tilde{g}=h$. By
 Cauchy-Schwartz inequality and (3.4), we have
 \begin{equation*} \tag{3.9}
 \begin{split}
 (\frac{\partial}{\partial
 t}-g^{\gamma\delta}\tilde{\nabla}_{\gamma}\tilde{\nabla}_{\delta})|g-\bar{g}|^{2}
 \leqslant&-2g^{\xi\eta}\tilde{g}^{\alpha\gamma}\tilde{g}^{\beta\delta}
 (\tilde{\nabla}_{\xi}g_{\alpha\beta}-\tilde{\nabla}_{\xi}\bar{g}_{\alpha\beta})
 (\tilde{\nabla}_{\eta}g_{\gamma\delta}-\tilde{\nabla}_{\eta}\bar{g}_{\gamma\delta})\\
 &
 +\frac{C}{\sqrt{t}}|g-\bar{g}|^{2}+C|\tilde{\nabla}(g-\bar{g})||g-\bar{g}|\\ \leqslant&
\frac{C}{\sqrt{t}}|g-\bar{g}|^{2}
\end{split}
 \end{equation*}
on $N^{n}\times [0,T_5]$.

Let $\varphi_1$  be the nonnegative function in Lemma 2.2 with
$a=1$, then
  \begin{eqnarray*}
 \frac{1}{C}(1+\tilde{d}(y,p))&\leqslant& \varphi_1(y)\leqslant C_0\tilde{d}(y,p)\ \ \  \text{on} \ N^{n}\backslash B(P,2),\\
 |\tilde{\nabla}\varphi_1|&+&|\tilde{\nabla}^{2}\varphi_1|\leqslant C,  \ \ \ \ \text{on} \
 N^{n}.
 \end{eqnarray*}
 For any fixed $t$ and any $\varepsilon>0$, consider the maximum of $|g-\bar{g}|^{2}-\varepsilon
 \varphi$. Clearly, the maximum is achieved at some
 point $P_{\varepsilon}^{t}$ and there hold
 \begin{eqnarray*}
|g-\bar{g}|^{2}(P_{\varepsilon}^{t})&\geqslant&
|g-\bar{g}|^{2}(y)-\varepsilon\varphi(y),\\
|\tilde{\nabla}|g-\bar{g}|^{2}|(P_{\varepsilon}^{t})\leqslant
C\varepsilon
&,&\tilde{\nabla}_{\alpha}\tilde{\nabla}_{\beta}|g-\bar{g}|^{2}
(P_{\varepsilon}^{t})\leqslant
C\varepsilon \tilde{g}_{\alpha\beta}(P_{\varepsilon}^{t})
 \end{eqnarray*}
 for all $y\in N^{n}$. This gives
 \begin{equation*} \tag{3.10}
 \begin{split}
  \limsup_{\varepsilon\rightharpoonup0}|g-\bar{g}|^{2}
  (P_{\varepsilon}^{t})&=\sup|g-\bar{g}|^{2}\\
  g^{\alpha\beta}\tilde{\nabla}_{\alpha}\tilde{\nabla}_
  {\beta}|g-\bar{g}|^{2}(P_{\varepsilon}^{t})&\leqslant
C\varepsilon
\end{split}
 \end{equation*}
 by the equivalence of $g$ and $\tilde{g}$.

Define a function
 $$
|g-\bar{g}|^{2}_{max}(t)=\sup_{y\in N^{n}}|g-\bar{g}|^{2}(y,t).
 $$
 By (3.9), and (3.10), we have
\begin{eqnarray*}
\frac{d^{+}}{d t}|g-\bar{g}|^{2}_{max}(t)&\leqslant&
\frac{C}{\sqrt{t}}|g-\bar{g}|^{2}_{max}(t)
\end{eqnarray*}
and then
$$
|g-\bar{g}|^{2}_{max}(t)\leqslant
e^{C\sqrt{T}}|g-\bar{g}|^{2}_{max}(0)=0
$$
Therefore the proof of the Proposition 3.3 is completed.

$\hfill\#$

\subsection{Proof of the main theorem}
 Let
$g_{ij}(x,t)$ and $\bar{g}_{ij}(x,t)$ be two solutions to the
Ricci flow (1.1) with bounded curvature and with the same initial
data. We solve the corresponding harmonic map flow (3.5) and (3.6)
with the same target
$(N^{n},h_{\alpha\beta})=(M^{n},g_{ij}(\cdot,T))$ respectively. We
obtain two solutions $F(x,t)$ and $\bar{F}(x,t)$ which are
diffeomorphisms for $t\in[0,T_5]$, where $T_5>0$ depends only on
$n, k_0, T$. Then ${(F^{-1})}^{*}g$ and
${(\bar{F}^{-1})}^{*}\bar{g}$ are two solutions to the Ricci-De
Turck flow with the same initial value. It follows from
Proposition 3.3 that
 $${(F^{-1})}^{*}g={(\bar{F}^{-1})}^{*}\bar{g},$$
 on $N^{n}\times[0,T_5]$.  So in order
to prove $g_{ij}(x,t)\equiv\bar{g}_{ij}(x,t)$, we only need to
show $F\equiv\bar{F}$. Let
\begin{eqnarray*}
 V^{\alpha}(y,t)&=&g^{\beta\gamma}(\tilde{\Gamma}^{\alpha}_
 {\beta\gamma}-{\Gamma}^{\alpha}_{\beta\gamma})=-(\triangle
 F\circ F^{-1})^{\alpha}\\
 \bar{V}^{\alpha}(y,t)&=&\bar{g}^{\beta\gamma}(\tilde{\Gamma}
 ^{\alpha}_{\beta\gamma}-\bar{\Gamma}^{\alpha}_{\beta\gamma})=-
 (\bar{\triangle}\bar{F}\circ
 \bar{F}^{-1})^{\alpha}.
 \end{eqnarray*}
 be two one-parameter family of vector fields on $N^{n}$,
 where ${g}_{\alpha\beta}(y,t)=({({F}^{-1})}^{*}{g})_
 {\alpha\beta}(y,t)$ and
  $\bar{g}_{\alpha\beta}(y,t)=({(\bar{F}^{-1})}^{*}\bar{g})_
  {\alpha\beta}(y,t)$. By
 Proposition 3.3, we have $g_{\alpha\beta}(y,t)=\bar{g}_
 {\alpha\beta}(y,t)$, thus the
 vector fields $V\equiv \bar{V}$ on the target $N^{n}$.
 Therefore, $F$ and $\bar{F}$ satisfy
 the same ODE equation with the same initial value:
  \begin{eqnarray*}
\frac{\partial}{\partial t}F&=&V\circ F,\\
 F(\cdot,0)&=&identity,
\end{eqnarray*} and
 \begin{eqnarray*}
 \frac{\partial}{\partial t}\bar{F}&=&V\circ \bar{F},\\
 \bar{F}(\cdot,0)&=&identity,
 \end{eqnarray*}
 By the same calculation as in the proof of Proposition 3.2, we have
 \begin{eqnarray*}\frac{\partial}{\partial t}{d}_{N^{n}}(F(x,t),\bar{F}(x,t))
 &\leqslant& \sup_{y\in N^{n}}|\tilde{\nabla}V|(y,t)
 {d}_{N^{n}}(F(x,t),\tilde{F}(x,t))\\
&\leqslant& \frac{C}{\sqrt{t}}
 {d}_{N^{n}}(F(x,t),\tilde{F}(x,t)).
 \end{eqnarray*}
This gives $${d}_{N^{n}}(F(x,t),\bar{F}(x,t))\leqslant
e^{C\sqrt{T}} {d}_{N^{n}}(F(x,0),\bar{F}(x,0))=0,$$
 which concludes that
 $$
 F(x,t)\equiv\bar{F(x,t)}.
 $$
 Thus  $g(x,t)=\bar{g}(x,t)$, for all $x\in M^{n}$ and
 $t\in [0,T_{5}]$ and for some $T_{1}>0$. Clearly, we can extend the
 interval $[0,T_{1}]$ to the whole $[0,T]$ by continuity method.

  Therefore we complete the proof of the Theorem 1.1.

  $\hfill\#$

  Finally, Corollary 1.2 is a direct consequence of Theorem 1.1.
  Indeed, since $G$ is the isometry group of $g_{ij}(x,0)$, then for any
  $\sigma\in G$, $\sigma^{*}g(\cdot,t)$ is still a solution to the
  Ricci flow with bounded curvature and $\sigma^{*}g(\cdot,t)\mid_{t=0}=\sigma^{*}g(\cdot,0)=g(\cdot,0).$
   By applying Theorem 1.1, we have
   $\sigma^{*}g(\cdot,t)=g(\cdot,t)$, $\forall t\in[0,T]$. So the
   corollary follows.

   $\hfill\#$

\end{document}